\documentclass{amsart}

\usepackage{amssymb}
\usepackage{amsfonts}
\usepackage{amsmath}
\usepackage{lastpage}
\usepackage[utf8]{inputenc}
\usepackage[legalpaper,bookmarks=true,colorlinks=true,linkcolor=blue,citecolor=blue]{hyperref}
\usepackage{graphicx}%
\usepackage{fancyhdr}
\usepackage{color}
\usepackage[mathlines]{lineno}
\usepackage{lscape}
\usepackage{epsfig}
\usepackage{natbib}
\usepackage{geometry}
\usepackage{tgbonum}
\fontfamily{qcr}\selectfont

\newtheorem{theorem}{Theorem}
\theoremstyle{plain}

\newtheorem{def-theo}{Theorem-Definition}
\newtheorem{def-pro}{Proposition and Definition}

\newtheorem{proposition}{Proposition}

\numberwithin{equation}{section}

\newcommand{\Bin}{\bigskip \noindent}

\newcommand{\Ni}{\noindent}

\begin{document}

\title{Second order Expansions for Extreme Quantiles of Burr Distributions and Asymptotic Theory of Record Values}

\author{Moumouni Diallo $^{(1)}$} 
\author{Modou Ngom $^{(2)}$}
\author{Soumaila Dembele $^{(3)}$}
\author{Gane Samb Lo $^{(4)}$}

\begin{abstract}
In this paper we investigate the Burr distributions family which contains twelve members. Second order expansions of quantiles of the Burr's distributions are provided on which may be based statistical methods, in particular in extreme value theory. Beyond the proper interest of these expansions, we apply them to characterize the asymptotic laws of their records of Burr's distributions, lead to new statistical tests.\\

\Ni $^{(1)}$ Moumouni Diallo.\\
Main Affiliation: Universit\'e des Sciences Sociales et Sciences de Gestion de Bamako (USSG-B), Mali.\\
Affiliation : LERSTAD, Gaston Berger University, Saint-Louis, S\'en\'egal\\
Imhotep Mathematical Center (IMC), imhotepsciences.org,\\
Email :moudiallo1@gmail.com\\

\Ni $^{(2)}$ Modou Ngom.\\
LERSTAD, Gaston Berger University, Saint-Louis, S\'en\'egal.\\
Imhotep Mathematical Center (IMC), imhotepsciences.org\\
Email: ngom.modou1@ugb.edu.sn\\

\Ni $^{(3)}$ Soumaila Dembele.\\
Universit\'e des Sciences Sociales et Sciences de Gestion de Bamako (USSG-B), Mali.\\
Affiliation : LERSTAD, Gaston Berger University, Saint-Louis, S\'en\'egal\\
Imhotep Mathematical Center (IMC), imhotepsciences.org\\
Emails : soumailadembeleussgb@gmail.com, dembele.soumaila@ugb.edu.sn\\

\Ni $^{(4)}$ Gane Samb Lo.\\
LERSTAD, Gaston Berger University, Saint-Louis, S\'en\'egal (main affiliation),\\
Department of Pure and Applied Mathematics, African university of Science and Technology, Abuja, Nigeria,\\
LSTA, Pierre and Marie Curie University, Paris VI, France (associated researcher),\\
Imhotep Mathematical Center (IMC), imhotepsciences.org,\\
Emails: gane-samb.lo@ugb.edu.sn , gslo@aust.edu.ng, ganesamblo@ganesamblo.net.\\

\Ni \textbf{keywords and phrases}: Burr distribution, Quantile expansion, record values, records times, extreme value theory, Gaussian asymptotic laws.\\
\textbf{AMS classification codes}: 62G30;60Fxx.
\end{abstract}
\maketitle

\newpage
\section{Introduction} \label{recEvtBeyond_01}

\Bin This papers deals of asymptotic laws of records values, extreme value theory when applied to the important family of Burr distributions (\cite{burr42}, \cite{burr68} and \cite{burr73}). Let us introduce to each of these three elements of the paper.\\

\Ni \textbf{Records theory}. The notion of records is present in real life at any corner. It is said that the year 2021 is the hottest one in History. In general, Records of many natural phenomena are monitored on a regularly basis: the coldest or hottest day, month, year, etc,; the rainiest month, year, etc.; the month or year with the greatest or smallest number of car or planes crashes, of biggest or smallest gross domestic product (GDP) [four countries]. At least any superlative corresponds to a record (lower record for positive superlatives, upper record for naive ones). Generally, records are associated to catastrophes or to big
 successes. Sports is punctuated by beaten records, in Olympic games, World championships. For example, the (lower records) of 100m in Athletics is particularly followed. Fr upper records in Sports, we can  cite the upper apnea record (time spent under water). So, the notion of record is extremely present in real file and it's modeling is highly valuable. 

\Ni \textbf{Univariate Extreme value Theory (UEVT)}. That theory is strongly connected to Records theory. Given a series of data $(X_j)_{j \geq 1}$, the \textit{UEVT} mainly studies the behavior of the partial maxima $M_n=max(X_1,\cdots,X_n)$, $n\geq 1$. Of course, each $M_{n+1}$ is a strong upper record value if and only if it exceeds the preceding one, that is $M_{n+1}>M_{n}$. The importance of \textit{UEVT} resides in the following paradox. It happens that some extreme events, for example
$p=\mathbb{P}(X>x)$ for large values of $x$, are not observed in samples and so, any plug-in estimator is exactly zero. In that context, how can we estimate the probability of occurrence of such events. The probability is usually given in the form $1/T$, where $T$ is defined as the temps needed to see a new occurrence of the event of exceedance of $x$. For large values of $x$, T is counted in thousands or more. In conclusion, the \textit{UEVT} is the theory of rare events. It's applications are countless and very important to circumvent catastrophes.\\

\Ni As expected, the asymptotic law of records values is strongly influenced by the asymptotic behavior of the extreme value. On that basis, we will give a more detailed account for these two theories but still concise in Subsections \ref{subUEVT} and \ref{subRecord}, respectively.\\

\Ni These two asymptotic theories are applied to the family of Burr's statistics whose elements have very interesting statistical properties and have important applications in a significant number of disciplines as in Telecommunications, Reliability, Actuarial Sciences,  Survival analysis, etc. So, we have to introduce to that family in Subsection \ref{subBurr}.\\


\subsection{Univariate Extreme value theory} \label{subUEVT} Let $X$, $X_1$, $X_2$, $\cdots$ be a sequence of independent real-valued randoms, defined on the same probability space $(\Omega, \mathcal{A},\mathbb{P})$, with common cumulative distribution function $F$, which has the lower and upper endpoints, the first asymptotic moment function and the generalized inverse function respectively defined by

$$
lep(F)=\inf \{x \in \mathbb{R}, \ F(x)>0\}, \  uep(F)=\sup \{x \in \mathbb{R}, \ F(x)<1\} 
$$

\noindent 

$$
R(x,F)=\frac{1}{1-F(x)}\int_{x}^{uep(F)}(1-F(y)) \ dy, \ x \in ]lep(F),uep(F)[
$$

\noindent and

$$
F^{-1}(u)=\inf \{x \in \mathbb{R}, \ F(x)\geq u\} \ for \ u \in ]0,1[ \ and \ F^{-1}(0)=F^{-1}(0+).
$$

\Bin $F$ is said to be in the extreme value domain of attraction of a non-degenerate \textit{df} $M$ whenever there
exist real and nonrandom sequences \ $(a_{n}>0)_{n\geq 1}$ and $(b_{n})_{n\geq 1}$ such that for any continuity point $x$\ of $M,$%

\begin{equation}
\lim_{n\rightarrow \infty} \mathbb{P}\left(\frac{X_{n,n}-b_{n}}{a_{n}}\leq
x\right)=\lim_{n\rightarrow \infty }F^{n}(a_{n}x+b_{n})=M(x).  \label{dl05}
\end{equation}

\Ni It is known that $M$ is necessarily of the family of the Generalized Extreme
Value (GEV) \textit{df} : 

\begin{equation}
H_{\gamma }(x)=\exp (-(1+\gamma x)^{-1/\gamma })\text{, }1+\gamma x\geq 0, \label{GEV}
\end{equation}

\Bin parametrized by $\gamma \in \mathbb{R}$, with $H_{0}(x)=1-\exp(-e^{-x})$, $x \in \mathbb{R}$, for $\gamma=0$. The parameter $\gamma $ is called
the extreme value index.\\

\Ni In this paper, we use some important facts from \textit{EVT} that we can summarize below, especially regarding functional representation of \textit{cdf}'s and their quantile functions in the extreme domain of attraction as well as their quantile functions.More details can be found in \cite{lo2021ALR}, as  a quick introduction on \textit{EVT} and to have  abroad view on how to find the domain of attraction of a \textit{cdf}. We will need the two following  propositions.\\

\begin{proposition} \label{proplo86} (see \cite{lo86}) We have the following equivalences.\\

\Bin Let $G(x)=F(e^{x}) $ the $cdf$ of the log-transformation. We have :\\

\noindent (1) If $\gamma>0$,

$$
F \in D(H_{\gamma}) \ \Leftrightarrow\ (G \in D(H_{0}) \ and \ R(x,G)\rightarrow \gamma \ as \ x \rightarrow uep(G)).
$$

\noindent (2) If $\gamma=0$,

$$
F \in D(H_{0}) \ \Leftrightarrow\ (G \in D(H_{0}) \ and \ R(x,G)\rightarrow 0 \ as \ x \rightarrow uep(G)).
$$

\noindent (3) If $\gamma<0$,

$$
F \in D(H_{\gamma}) \ \Leftrightarrow\ G \in D(H_{\gamma}). 
$$

\end{proposition}

\Ni The next proposition provides interesting functional representations of quantile functions of $cdf$'s in the extreme domain of attraction.

\begin{proposition} \label{portal.rd} (\cite{karamata30} and \cite{dehaan}) We have the following characterizations for the three extremal domains.

\bigskip \noindent (a) $F\in D(H_{\gamma })$, $\gamma >0$, if and only if there exist a
constant $c$ and functions $a(u)$ and $\ell (u)$ of  $u\in
]0,1[$ satisfying

\begin{equation*}
(a(u),\ell (u))\rightarrow (0,0)\text{ as }u\rightarrow 0,
\end{equation*}

\bigskip \noindent such that $F^{-1}$ admits the following representation of Karamata

\begin{equation}
F^{-1}(1-u)=c(1+a(u))u^{-\gamma }\exp \left(\int_{u}^{1}\frac{\ell (t)}{t}dt\right). \label{portal.rdf}
\end{equation}

\bigskip \noindent (b) $F\in D(H_{\gamma })$, $\gamma <0,$ if and only if $uep(F)<+\infty $ and
there exist a constant $c$ and functions $a(u)$ and $\ell (u)$ of $u\in ]0,1[$ satisfying
\begin{equation*}
(a(u),\ell (u))\rightarrow (0,0)\text{ as }u\rightarrow 0,
\end{equation*}

\bigskip \noindent such that $F^{-1}$ admit the following representation of Karamata

\begin{equation}
uep(F)-F^{-1}(1-u)=c(1+a(u))u^{-\gamma }\exp \left(\int_{u}^{1}\frac{\ell (t)}{t}
dt\right). \label{portal.rdw}
\end{equation}

\bigskip \noindent (c) $F\in D(H_{0})$ if and only if there exist a constant $d$ and a slowly
varying function $s(u)$ such that 

\begin{equation}
F^{-1}(1-u)=d+s(u)+\int_{u}^{1}\frac{s(t)}{t}dt,0<u<1, \label{portal.rdg}
\end{equation}

\bigskip \noindent and there exist a constant $c$ and functions $a(u)$ and $\ell (u)$ of $u\in ]0,1]$ satisfying
\begin{equation*}
(a(u),\ell (u))\rightarrow (0,0)\text{ as }u\rightarrow 0,
\end{equation*}

\bigskip \noindent such that the function $s(u)$ of $u \in ]0,1[$ admits the representation
\begin{equation}
s(u)=c(1+a(u)) \exp \left(\int_{u}^{1}\frac{\ell (t)}{t}dt\right). \label{portal.rdgs}
\end{equation}

\noindent Moreover, if $F^{-1}(1-u)$ is differentiable for small values of $u$ such
that $r(u)=-u(F^{-1}(1-u))^{\prime }= -u\ dF^{-1}(1-u)/du$ is slowly varying at
zero, then (\ref{portal.rdg}) may be replaced by 

\begin{equation}
F^{-1}(1-u)=d+\int_{u}^{u_{0}}\frac{r(t)}{t}dt,0<u<u_{0}<1, \label{portal.rdgr}
\end{equation}

\noindent which will be called a \textit{reduced de Haan representation} of $F^{-1}.$
\end{proposition}

\Bin On top of these representations, we may use the simple criteria (See more criteria in \cite{dehaan} and \cite{lo86}, recalled in theorem \ref{theoALL1}).

\begin{proposition} \label{criteria} We have:\\

\Ni (a)  $F \in D(G_{\gamma})$, $\gamma>0$ if and only if \textit{uep}$(F)=+\infty$ and

$$
\forall \lambda>0, \; u \in ]0,1[, \ \lim_{u \rightarrow 0} \frac{F^{-1}(1-\lambda u)}{F^{-1}(1- u)}=\lambda^{-\gamma}
$$

\Bin (b)  $F \in D(G_{\gamma})$, $\gamma<0$ if and only if $\textit{\text{uep}}(F)<\infty$ and

$$
\forall \lambda>0, u \in ]0,1[, \; \lim_{u \rightarrow 0} \frac{uep(F)-F^{-1}(1-\lambda u)}{uep(F)-F^{-1}(1- u)}=\lambda^{-\gamma}
$$

\Bin (c)  $F \in D(G_{\gamma})$, $\gamma=0$ if and only there exists a function $s(u)$ of $u\in ]0,1[$ which is slowly varying function at zero and such that 

$$
\forall \lambda>0, u \in ]0,1[, \ \lim_{u \rightarrow 0} \frac{F^{-1}(1-\lambda u)-F^{-1}(1-u)}{s(u)}=-\log \lambda
$$

\Bin if and only if

$$
\forall \lambda>0, \forall  \;  0<\mu\neq 1, u \in ]0,1[, \ \lim_{u \rightarrow 0} \frac{F^{-1}(1-\lambda u)-F^{-1}(1-u)}{F^{-1}(1-\mu u)-F^{-1}(1-u)}=
\frac{\log \lambda}{\log \mu}.
$$
\end{proposition}

\Bin
\subsection{Gaussian asymptotic laws of record values}\label{subRecord}

\Ni Let us begin by define records values and records times for a sequence of real random variables defined on the same probability space $\left(\Omega, \mathcal{A}, \mathbb{P}\right)$: $Y_1$, $Y_2$, $\cdots$. Record times and record values are defined as follows.\\

\Ni \textbf{Strong upper record times}. Let us put $u(1)=1$ as the first strong upper record time. For any $n\geq 2$, we define, by induction, whenever the $(n-1)^{th}$ upper record time $u(n-1)$ exists,

$$
U_{n}=\left\{j>u(n-1), \ Y_j>Y_{u(n-1)}\right\}.
$$

\Bin Hence, for $n\geq 2$, the $(n)^{th}$ upper record time is $u(n)=+\infty$ if $U_n$ is empty and, otherwise

$$
u(n)=\inf U_n.
$$

\Bin  \textbf{Strong lower record times}. Let us put $\ell(1)=1$ as the first strong lower record time. For any $n\geq 2$, we define, by induction, whenever the $(n-1)^{th}$ lower record time $\ell(n-1)$ exists,

$$
L_{n}=\left\{j>\ell(n-1), \ Y_j<Y_{\ell(n-1)}\right\}
$$

\Bin  Hence, for $n\geq 2$, the $n$-\textit{th} lower record time is $\ell(n)=+\infty$ if $L_n$ is empty and, otherwise

$$
\ell(n)=\inf L_n.
$$

\Bin  \textbf{Strong record values}. For each $n\geq 1$ such that $u(n)$ is finite, we have a sequence of  strong upper record values  
$$(Y^{(k)}=Y_{u(k)}, \ 1\leq k\leq n).$$

\Ni For each $n\geq 1$ such that $\ell(n)$ is finite, we have a sequence of  strong lower record values  
$$Y_{(k)}=Y_{\ell(k)}, \ 1\leq k\leq n.$$

\Bin \textit{There are many results on probability laws of record values and record times, especially for  \textit{iid} random variables with common 
{\textit{cdf} $F$, eventually associated with the \textit{pdf} $f$ with respect to the Lebesgue measure $\lambda$}[ or \textit{iid} random variables with common mass probability functions $p$]}. Important books introduce the the study of records in the \textit{iid} scheme, such as \cite{nevzorov}, \cite{ahsan88}, \cite{ahsan95}, \cite{ahsan04}, \cite{ahsanullahUGB}, \cite{arnold}, etc. Also, statistical applications are largely available and  Characterization problems involving functional equations (See \cite{ahsan88}, \cite{nevzorov}, \cite{resnick}, etc.).\\

\Ni In a recent paper, we are interesting in finding the asymptotic laws of records values of elements of Burr's family. In that extent, we will mainly follow  \cite{lo2021ALR} who provided practical methods of finding the asymptotic law of record values depending, in general, on the extreme value attraction domain of a distribution $F$. We intend to use such results to the Burr's family. Let us broaden the notation in order to be able to expose the main theorem of the cited authors. Given the sequence defined above, we consider the sequence of strong record values $X^{(1)}=X_1$, $X^{(n)}$, $\cdots$ and the sequence of record times  $U(1)=1$, $U(2)$, $\cdots$. Some conditions depend on a sequence of  finite sum of $n$ standard exponential random variables, $n\geq 1$,

$$
S_{(n)}=E_{1,n} + \cdots + E_{n,n},
$$

\noindent and we denote

$$
V_n=\exp(-S_{(n)}) \ and \ v_n=\exp(-n), \ n\geq 1.
$$

\noindent From this, we may set

$$
(Ha) \ : \ \sup\left\{ \left|\frac{s(u)}{s(v)} - 1 \right|, \ \min(v_n,V_n) \leq u,v \leq \max(v_n,V_n)\right\} \rightarrow_{\mathbb{P}} 0 \ as \ n\rightarrow +\infty,
$$

\noindent

$$
(Hb) \ : \ (\exists \alpha>0), \ \sqrt{n} \ s(v_n) \rightarrow \alpha \ as \ n\rightarrow +\infty,
$$

\noindent \noindent where $\rightarrow_{\mathbb{P}}$ stands for the convergence in probability.\\

\noindent \cite{lo2021ALR} obtained the results below  that cover the whole extreme value domain of attraction. Let us begin by asymptotic laws for $F \in \mathcal{D}=D({G_{\gamma}})$, $\gamma \in \mathbb{R}$.

\begin{theorem} \label{theoALL1} Let $F \in \mathcal{D}$. We have :\\

\noindent (a) If $\gamma>0$, the asymptotic law of $X^{(n)}$ is lognormal, precisely

$$
\left(\frac{X^{(n)}}{F^{-1}\left(1-e^{-n}\right)}\right)^{n^{-1/2}} \rightsquigarrow LN(0,\gamma^2),
$$

\noindent where $LN(m,\sigma^2)$ is the lognormal law of parameters $m$ and $\sigma>0$.\\

\noindent (b) If $\gamma>0$ and $X > 0$, $Y=\log X \in D(G_{\gamma})$ and $R(x,G)\rightarrow \gamma \ as \ x \rightarrow uep(G)$ and we have

$$
\frac{Y^{(n)}- G^{-1}\left(1-e^{-n}\right)}{\sqrt{n}} \rightsquigarrow \mathcal{N}(0,\gamma^2).
$$

\noindent (c) If $\gamma<0$, the asymptotic law of $X^{(n)}$ is lognormal, precisely

$$
\left(\frac{uep(F)-X^{(n)}}{uep(F)-F^{-1}\left(1-e^{-n}\right)}\right)^{n^{-1/2}} \rightsquigarrow \exp(\mathcal{N}(0,\gamma^2)).
$$

\bigskip \noindent (d) Suppose that $\gamma=0$ and  $R(x,G)\rightarrow 0 \ as \ x \rightarrow uep(G)$. If (Ha) and (Hb) hold both, we have

$$
X^{(n)}- F^{-1}\left(1-e^{-n}\right) \rightsquigarrow \mathcal{N}(0,\alpha^2).
$$

\noindent More precisely, we have :  Given $\gamma=0$,  $R(x,G)\rightarrow 0 \ as \ x \rightarrow uep(G)$ and (Ha), the above asymptotic normality is valid if and only if (Hb) holds.
\end{theorem}

\noindent The theorem can be extended outside the extreme domain of attraction as follows. Suppose that:\\

\noindent (Ga) $uep(F)=+\infty $ and $F$ is differentiable in some neighborhood of $]x_0, \ +\infty[$.\\

\noindent (Gb) The function 

$$
s(x)=e^{-x} F^{-1}(1-t)\biggr|_{t=e^{-x}}, \ e^{x}<u_0<1, \ for \ some \ u_0 \in ]0,1[
$$ 

\noindent \textbf{decreases} to $0$ as $x \rightarrow +\infty$ and is : for any sequence $(x_n,y_n)_{n\geq 1}$ such that
$$
\limsup_{n\rightarrow +\infty} |x_n-y_n|/\sqrt{n} <+\infty,
$$

\noindent we have, for some $\alpha>0$,

$$
\lim_{n\rightarrow +\infty} \sqrt{n} \ s\left(\exp(\min(x_n,y_n))\right)=\lim_{n\rightarrow +\infty} \sqrt{n} \ s\left(\exp(\max(x_n,y_n))\right)=\alpha.
$$

\noindent We have the following generalization.

\begin{theorem} \label{theoALL2} If $F$ satisfies Assumptions (Ga) and (Gb), we have

$$
X^{(n)}- F^{-1}\left(1-e^{-n}\right) \rightsquigarrow \mathcal{N}(0,\alpha^2) .
$$
\end{theorem}

\Bin \textbf{Important result}. For $\gamma\neq 0$, we need no condition for the asymptotic law to hold true.\\

\bigskip \noindent \textbf{Rule of working \label{commentsI}}. Suppose that $F$ lies in $\mathcal{D}$, the extreme domain of attraction.\\

\noindent (e) If $F \in D(H_{\gamma})$, $\gamma\neq 0$, we apply Points (a) or (c) of theorem \ref{theoALL1} without any further condition.\\

\noindent (f) If $F \in D(H_{0})$ and $\exp(X) \in D(H_{\gamma})$ for some $\gamma>0$, we apply Point (b) without any further condition.\\

\noindent (g) If $F \in D(H_{0})$ and $s(u) \rightarrow 0$ as $u\rightarrow 0$ and if (Ha) and (Hb) hold, we conclude by applying Point (d). If not (as it is for a lognormal law), we search whether $X_1=\exp(X) \in D(H_{\gamma})$ for some $\gamma>0$ or $X_1=\exp(X)$ fulfills (Ha) and (Hb). If yes, we conclude by Point (b) or Point (d). If not, we consider $X_2=\exp(X_1)$, and we continue the process until we reach $X_p=\exp(X_{p-1}) \in D(G_{\gamma})$ for some $\gamma>0$ or $X_p=\exp(X_{p-1})$ for some $p\geq 1$.\\

\Ni Other results in \cite{lo2021ALR} are expressed as representations as in the following theorem.

\begin{theorem} \label{theoALL3} Let $F \in D(H_{\gamma})$, $\gamma \in \mathbb{R}$. Then, there exists a probability space $(\Omega, \mathcal{A}, \mathbb{P})$ holding a sequence of independent standard exponential random variables $(E_n)_{n\geq 1}$ and a Brownian Process $\{W(t), \ t\geq 0\}$ such that the record values $X^{(n)}$, $n\geq 1$,  of the sequence $X_j=F^{-1}\left(1-e^{E_j}\right)$, $j\geq 1$, satisfy the following representations below under the appropriate conditions. Here, $S_n=E_1+...+E_n$, $n\geq 1$, are the partial sums of the sequence $(E_n)_{n\geq 1}$, $S_n^{\ast}=n^{-1/2}(S_n-n)$,  $v_n=e^{-n}$ and $V_n=e^{-S_n}$. Below, the function $a(u)$, $b(u)$ and $s(u)$ of $u \in ]0,1[$ are those in the representations in Proposition \ref{portal.rd}.\\

\noindent By denoting $W_n^{\ast}=n^{-1/2}W(n)$ and $c_n=n^{-1/2} \log n$, we have

$$
W_n^{\ast} \sim \mathcal{N}(0,1) \ \ and \ \ \left|S_n^{\ast}-W_n^{\ast}\right|=O_{\mathbb{P}}(c_n).
$$

\bigskip \noindent Further, we have the following results.\\

\noindent (a) Let $\gamma>0$. Suppose that

\begin{equation}
1-\frac{1+a(V_n)}{1+a(v_n)}=O_{\mathbb{P}}(a_n), \ \sup\{|b(t)|, \ 0\leq t \leq v_n \vee V_n\}=O_{\mathbb{P}}(b_n) \ . \label{vitConv}
\end{equation}

\noindent Then, we have

\begin{eqnarray*}
\left(\frac{X^{(n)}}{F^{-1}\left(1-e^{-n}\right)}\right)^{n^{-1/2}}&=& \exp(\gamma S_n^{\ast}) + O_{\mathbb{P}}(a_n \vee b_n)\\
&=&\exp(\gamma W_n^{\ast})+ O_{\mathbb{P}}(a_n \vee b_n \vee c_n).
\end{eqnarray*}

\bigskip
\noindent (b) Let $\gamma>0$ and $X > 0$, $Y=\log X \in D(G_0)$ and $R(x,G)\rightarrow \gamma \ as \ x \rightarrow uep(G)$ and we have

\begin{eqnarray*}
\frac{Y^{(n)}- G^{-1}\left(1-e^{-n}\right)}{\sqrt{n}}&=&\gamma S_n^{\ast} + O_{\mathbb{P}}(a_n \vee b_n) \\
&=&\gamma W_n^{\ast}+ O_{\mathbb{P}}(a_n \vee b_n \vee c_n).
\end{eqnarray*}

\noindent (c) Let $\gamma<0$. Then, by using the rates of convergence in Formula \eqref{vitConv}, we have 

\begin{eqnarray*}
\left(\frac{uep(F)-X^{(n)}}{uep(F)-F^{-1}\left(1-e^{-n}\right)}\right)^{n^{-1/2}}&=& \exp(\gamma S_n^{\ast})+ O_{\mathbb{P}}(a_n \vee b_n )\\
&=& \exp(\gamma W_n^{\ast})+ O_{\mathbb{P}}(a_n \vee b_n \vee c_n).
\end{eqnarray*}

\bigskip \noindent (d) Suppose that $\gamma=0$ and  $R(x,G)\rightarrow 0 \ as \ x \rightarrow uep(G)$. Suppose that (Ha) and (Hb) hold both. If

$$
\sup\left\{ \left|\frac{s(u)}{s(v)} - 1 \right|, \ \min(v_n,V_n) \leq u,v \leq \max(v_n,V_n)\right\}=O_{\mathbb{P}}(d_n), \ and \ \sqrt{n} s(v_n)-\alpha=O(e_n),
$$

\bigskip
\noindent we have

\begin{eqnarray*}
X^{(n)}- F^{-1}\left(1-e^{-n}\right)&=&\alpha S_n^{\ast}+O_{\mathbb{P}}(d_n \vee e_n)\\ 
&=&\alpha W_n^{\ast} +O_{\mathbb{P}}(c_n \vee d_n \vee e_n) .
\end{eqnarray*}
\end{theorem}

\bigskip \noindent \textbf{Rules of working \label{commentsII}}. In the domain of extremal attraction, most of the \textit{cdf}'s which are used in applications are differentiable in a left-neighborhood of the upper endpoint. In such a case, we may take $a\equiv0$ in Representation \eqref{portal.rdf} and \eqref{portal.rdw} in Proposition \ref{portal.rd}. By solving easy differential equations, we have the representation for

\begin{equation}
b(u)=-u\biggr(G^{-1}(1-u)\biggr)^{\prime }-\gamma, \ u\in ]0,1[ \ \ and \ \ a\equiv 0
\end{equation}

\noindent  for $\gamma>0$ and

\begin{equation}
b(u)= -\gamma -\frac{u}{F^{\prime}\biggr(F^{-1}(1-u)\biggr)\biggr(uep(F)-F^{-1}(1-u)\biggr)}, \ u\in ]0,1[.
\end{equation}

\bigskip
\noindent for $\gamma<0$, whenever we have  $b(u)\rightarrow 0$ as $u\rightarrow 0$. Consequently, the rate of convergence reduces to $O_{\mathbb{P}}(b_n \vee c_n)$.\\

\noindent For $\gamma=0$, Representation \eqref{portal.rdgr} in Proposition \ref{portal.rd} holds for

$$ 
s(u)=-u\biggr(F^{-1}(1-u)\biggr)^{\prime }, \ 0<u<1,
$$

\bigskip \noindent whenever it is slowly varying at zero and the rate of convergence $d_n$ becomes useless. In such case, the rate of convergence reduces $O_{\mathbb{P}}(d_n \vee c_n)$.\\

\noindent Furthermore, based on the limit $S_n/n \rightarrow 1$ as $n\rightarrow +\infty$, we get that, for any $\eta \in ]0,1[$ ,

\begin{equation}
\liminf_{n \rightarrow +\infty} \mathbb{P}\biggr(e^{- n/\eta} \leq e^{-S_n} \leq e^{-\eta n}\biggr)=1. \label{approS-n}
\end{equation}

\bigskip
\noindent So we may replace the rates of convergence \textit{$d_n$ and $b_n$} by $d_n(\eta)$ and $b_n(\eta)$ defined as follows, for $\eta \in ]0,1[$

\begin{equation}
\sup\{|b(t)|, \ 0\leq t \leq e^{-\eta n}\}=O(b_n(\eta)) \label{bprime-n}
\end{equation}

\bigskip
\noindent and 

\begin{equation}
\sup\left\{ \left|\frac{s(u)}{s(v)} - 1 \right|, \ e^{- n/\eta} \leq u,v \leq e^{- \eta n} \right\}=O(d_n(\eta)). \label{dprime-n}
\end{equation}


\subsection{Burr's Family} \label{subBurr} We may quote \cite{devroye}, "In a serie of papers, Burr (\cite{burr42}, \cite{burr68} and \cite{burr73}) has proposed a versatile family of densities". That class has twelve (12) statistical distributions given in Table \ref{tab01}, in which $k$, $c$ and $r$ are positive parameters and the support or domain of each element of the family are precised. We added the distribution (Xa) as a version of the distribution (X). 

\begin{table}
\centering
\begin{tabular}{c||c|c|c}
\hline\hline
Number & $F(x)$ & Domain & $F^{-1}(u)$ \\ \hline\hline
I & $F(x)=x$ & $(0,1)$ & $F^{-1}(u) =u $ \\ 
&  &  &  \\ 
II & $\biggr(1+e^{-x}\biggr)^{-r}$ & $\mathbb{R}$ & $\log \biggr(\frac{%
u^{1/r}}{1-u^{1/r}}\biggr) $ \\ 
&  &  &  \\ 
III & $\biggr(1 + x^{-k}\biggr)^{-r}$ & $\mathbb{R}_+$ & $\biggr( \frac {%
u^{1/r}}{1-u^{1/r}} \biggr) ^{1/k}$ \\ 
&  &  &  \\ 
IV & $\biggr(1+\left(\frac{c-x}{x}\right)^{1/c}\biggr)^{-r}$ & $(0,c)$ & $%
\frac{c}{1+\biggr( \frac {1-u^{1/r}}{u^{1/r}}\biggr) ^{c}}$ \\ 
&  &  &  \\ 
V & $\biggr(1+k e^{-\tan x}\biggr)^{-r}$ & $(-\pi/2, \ \pi/2)$ & $\arctan %
\biggr(\log \biggr( \frac{k}{\left( u^{-1/r}-1\right) }\biggr) \biggr) $ \\ 
&  &  &  \\ 
VI & $\biggr(1+ k e^{- \sinh x}\biggr)^{-r}$ & $\mathbb{R}$ & $\arg \sinh%
\biggr( \log \left( \frac{k}{u^{-1/r}-1}\right) \biggr)$ \\ 
&  &  &  \\ 
VII & $2^{-r}\biggr(1+ \tanh x \biggr)^{r}$ & $\mathbb{R}$ & $\arg \tanh %
\biggr( 2u^{1/r}-1\biggr)$ \\ 
&  &  &  \\ 
VIII & $\biggr(\frac{2}{\pi} \arctan (e^{x})\biggr)^{r}$ & $\mathbb{R}$ & $%
\log \biggr( \tan \biggr( \frac{\pi }{2}u^{1/r}\biggr) \biggr) $ \\ 
&  &  &  \\ 
IX & $1 - \biggr(\frac{2}{2+k\left(\left(1+e^x\right)^r-1\right)}\biggr)$ & $%
\mathbb{R}$ & $\log \biggr( \biggr( 1+k^{-1}\biggr( \frac{1+u}{1-u}\biggr) %
\biggr) ^{1/r}-1 \biggr) $ \\ 
&  &  &  \\ 
X & $\biggr(1+e^{-x^2}\biggr)^{-r}$ & $\mathbb{R}_+$ & $\biggr( \log \left( 
\frac{1}{u^{-1/r}-1}\right) \biggr) ^{1/2}$ \\ 
&  &  &  \\ 
XI & $\biggr(x - \frac{1}{2\pi} \sin(2\pi x)\biggr)^{r}$ & $(0,1)$ & non
explicit \\ 
&  &  &  \\ 
XII & $1 - \biggr(1+ x^c\biggr)^{-r}$ & $\mathbb{R}_+$ & $\biggr(\biggr( 1-u%
\biggr) ^{-1/r}-1\biggr) ^{1/c}$ \\ 
&  &  &  \\ 
Xa & $\biggr(1-e^{-x^2}\biggr)^{r}$ & $\mathbb{R}_+$ & $\biggr( \log \left( 
\frac{1}{1-u^{1/r}}\right) \biggr) ^{1/2}$ \\ \hline\hline
\end{tabular}
\caption{Burr's distributions}
\label{tab01}
\end{table}

\Bin Since, some elements of that family have taken important roles in many parts of Statistics and have been extended a great number of times. Let us denote any Burr distribution by $\mathcal{B}(T,a,b,c)$ where $T$ stands of (I), $\cdots$, (XII), the last cited parameter is the final power of the \textit{cdf}, if such a power exists, and the first parameters are cited in the order of appearance.\\

\Ni Also, that family intersects  with celebrated other families of distribution : Pearson, Dagum and Singh-Magdalla families to cite a few. For example the Sing-Madalla (See \cite{singhMad}, \cite{dagum}, \cite{rasheed}, etc ) defined by

$$
F_{sm(a,b,c)}(x)=1 -(1+ax^c)^{-r}, \ x\geq 0, \ \ a>0, \ b>0, \ c>0
$$ 

\Bin reduces to the $\mathcal{B}(\textit{\text{XII}},1,c,r)$ law. If $X \sim   sm(a,b,c) $ , the variable $1/X$ becomes a Dagum law with

$$
F_{D(a,b,c)}(x)=(1+ax^{-b})^{-c}, \ x\geq 0, \ \ a>0, \ b>0, \ c>0.
$$

\Bin That element of the Dagum class has been generalized in \cite{rasheed} as the \textit{Topp-Leone Dagum Distribution} of parameters $a>0 $, $b>0$, $c>0$, $d>0$, $f>0$,

\begin{equation}
F(x)=\biggr( 1 -\left( \left\{1-(1+ax^{-b})^{-c}\right\}\right)^{d}\biggr)^{f}, \ x\geq 0, \label{tlDagum}
\end{equation}

\Bin where $f=2$ is \cite{rasheed}. But we let $f>0$ and we still have a \textit{cdf}.\\

\Ni Especially, the $\mathcal{B}(\textit{\text{XII}},a,b,c)$ distribution is an instrumental tool in extreme value distribution (See for example \cite{Segers}). 
Also the $\mathcal{B}(\textit{\text{III}},a,b,c)$ distribution is also a member of the Dagum distributions system (\cite{kleiber}) which Dagum himself named as a \textit{generalized Burr system}.\\

\subsection{Motivations and organization of the paper}
\Bin Our achievements in that paper is that we entirely characterized the asymptotic laws for record values of \textit{cdf}'s in the Burr family based on the second order expansions of their quantile functions. Statistical tests are derived from these results. Simulations are given to back the results. But adaptive tests will not be considered here. They will be the object of an applied statistics paper.\\

\Ni The rest of the paper is organized as follows. In Section \ref{sec_02}, we expose the second order expansions of their quantile functions of \textit{cdf}'s in the Burr family the extremal domains of attraction and the asymptotic law of their record values. In Section \ref{sec_03}, we proceed to a simulation study. The very technical part on the second order expansions of quantile functions in given in Section \ref{sec_04}. A conclusion part (Section \ref{sec_05}) finishes the paper.

\section{Our main results} \label{sec_02}

For each \textit{cdf} element of the Burr family, we give the expansion of the quantile function, determine the extreme value domain $D(G_{\gamma})$ and next apply Theorems 1, 2 and 3 in \cite{lo2021ALR} to find the asymptotic law of the record values. However, we need to complete their theorem 3 by Theorem 
\ref{extensionBelow} below.\\

\begin{theorem} \label{extensionBelow}
Suppose that $\gamma=0$ and  $R(x,G)\rightarrow 0 \ as \ x \rightarrow uep(G)$. Suppose that (Ha) holds. If

$$
\sup\left\{ \left|\frac{s(u)}{s(v)} - 1 \right|, \ \min(v_n,V_n) \leq u,v \leq \max(v_n,V_n)\right\}=O_{\mathbb{P}}(d_n), 
$$

\bigskip \noindent then we have

\begin{eqnarray*}
\frac{X^{(n)}-F^{-1}(1-e^{-n})}{s(v_n)\sqrt{n}}=S_n^{\ast} + O_{\mathbb{P}}\left(d_n\right)=W_n^{\ast} + O_{\mathbb{P}}(c_n \vee d_n).
\end{eqnarray*}

\end{theorem}

\Ni \textbf{Proof}. In the proof of Theorem 3 in \cite{lo2021ALR}, in Formula (25), we have \\

\begin{eqnarray*}
X^{(n)}-F^{-1}(1-e^{-n})&=&s(V_n)-s(v_n)+\int_{v_n}^{V_n} \frac{s(u)}{u} \ du \label{pc01}\\
&=&s(v_n) \left(\frac{s(V_n)}{s(v_n)}-1\right)+s(v_n)(S_{n}-n)(1+O_{\mathbb{P}}(d_n)).\notag
\end{eqnarray*}

\Bin Hence

\begin{eqnarray*}
\frac{X^{(n)}-F^{-1}(1-e^{-n})}{s(v_n)\sqrt{n}}&=&\frac{1}{n^{1/2}} \left(\frac{s(V_n)}{s(v_n)}-1\right)+  (1+O_{\mathbb{P}}(d_n)) \ S_n^{\ast}\\
&=&S_n^{\ast} + \frac{O_{\mathbb{P}}(d_n)}{n^{1/2}} + S_n^{\ast} O_{\mathbb{P}}(d_n)=S_n^{\ast} + O_{\mathbb{P}}(d_n) .
\end{eqnarray*}

\Bin So the conclusion

\begin{eqnarray*}
\frac{X^{(n)}-F^{-1}(1-e^{-n})}{\sqrt{n}s(v_n)}=S_n^{\ast} + O_{\mathbb{P}}\left(d_n\right)=W_n^{\ast} + O_{\mathbb{P}}(c_n \vee d_n) ,
\end{eqnarray*}

\Bin is straightforward. Of course if $n^{1/2} s(v_n)\rightarrow \alpha>0$, we get the result of Theorem \ref{theoALL2} again. $\blacksquare$\\

\Ni Now, let us state the asymptotic laws of record values for all members of the Burr family and the distribution (Xa) in the theorem below. All the proofs of the  quantile function of second order expansions are given in Section \ref{sec_04}. For every asymptotic law of record values, we precise the arguments 
(from Theorems \ref{theoALL1} and/or \ref{theoALL3}, and/or \ref{extensionBelow}) to be applied. We do not need to give the detailed proof of all of them. However, after the statement of all the results, for each argument in Theorems \ref{theoALL1} and/or \ref{theoALL3}, and/or \ref{extensionBelow}, we will give the proof of a case on which it is applied.\\
 
\begin{theorem}\label{theoI}  Let $X$ follows $\mathcal{B}(T,a,b,c)$ with \textit{cdf} $F$. We have the following results for the distribution in the table \ref{tab01}.\\

\Ni \textbf{Burr I}.\\

\Ni (i) The quantile function is expanded as follows.\\ 
 
\begin{equation}
uep(F)-F^{-1}(1-u)=u. \label{ResquantDefI}
\end{equation}

\Bin (ii) $F \in D(G_{\gamma})$, $\gamma=-1$, $uep(F)=1$.\\

\Ni (iii) The asymptotic law of the record values $X^{(n)}$, $n\geq 1$, is given as follows.

\begin{eqnarray}
\biggr(e^{n}(1-X^{(n)})\biggr)^{n^{-1/2}}&=&\exp(-S_{n}^{\ast})=\exp(-W_{n}^{\ast})+ O_{\mathbb{P}}(c_n)\notag\\
	&\rightarrow& \exp(\mathcal{N}(0,1)). \label{al_records_Burr_I}
\end{eqnarray}

\Bin \textbf{Burr II}.\\

\Ni (i) The quantile function is expanded as follows.\\

\begin{equation}
F^{-1}(1-u)=\log r + \log (1/u) - \frac{r+1}{2r} u + O(u^2).
\label{ResquantDefII}
\end{equation}

\Bin (ii) Here $Z=\exp(X) \in D(G_{\gamma})$, $\gamma=1$, i.e., $X=\log Z$ and $uep(F)=+\infty$. 

\Bin (iii) The asymptotic law of the record values $X^{(n)}$, $n\geq 1$, is given as follows.

\begin{eqnarray}
\frac{X^{(n)} - n}{\sqrt{n}}&=&S_{n}^{\ast}+O_{\mathbb{P}}(d_n(\eta))= W_{n}^{\ast}+ O_{\mathbb{P}}(d_n(\eta)\vee c_n) \notag\\
	&\rightarrow& \mathcal{N}(0,1). \label{al_records_Burr_II}
\end{eqnarray}

\Bin \textbf{Burr III}.\\

\Ni (i) The quantile function is expanded as follows.\\ 
 
\begin{equation}
F^{-1}(1-u)= r^{1/k} u^{-1/k} \biggr(1 - \frac{r+1}{2kr} u + O(u^2)\biggr).
\label{ResquantDefIII}
\end{equation}

\Bin (ii) $F \in D(G_{\gamma})$, $\gamma=1/k$, $uep(F)=+\infty$.\\

\Ni (iii) The asymptotic law of the record values $X^{(n)}$, $n\geq 1$, is given as follows.

\begin{eqnarray}
\biggr(\frac{X^{(n)}}{r^{1/k} \exp(n/k)}\biggr)^{n^{-1/2}}&=&\exp\left(\frac {1}{k} S_{n}^{\ast}\right)+O_{\mathbb{P}}(a_n\vee b_n)\notag\\
&=& \exp\left(\frac {1}{k} W_{n}^{\ast}\right)+ O_{\mathbb{P}}(a_n\vee b_n \vee c_n)\notag\\
	&\rightarrow& \exp(\mathcal{N}(0,k^{-2})). \label{al_records_Burr_III}
\end{eqnarray}

\Bin \textbf{Burr IV}.\\

\Ni (i) The quantile function is expanded as follows.\\ 
 
\begin{equation}
c-F^{-1}(1-u)=cr^{-c}u^{c}\left( 1+\frac{c(r+1)}{2r}u+O\left(
u^{2}\right) \right). \label{ResquantDefIV}
\end{equation}

\Bin (ii) $F \in D(G_{\gamma})$, $\gamma=-c$, $uep(F)=c$.\\

\Ni (iii) The asymptotic law of the records value $X^{(n)}$, $n\geq 1$, is given as follows.

\begin{eqnarray}
\biggr( \frac{r^{c} e^{cn}}{c}  \left(c - X^{(n)}\right)\biggr)^{n^{-1/2}}&=&\exp(-c S_{n}^{\ast})+O_{\mathbb{P}}(a_n\vee b_n) \notag\\
&=& \exp(-c W_{n}^{\ast})+ O_{\mathbb{P}}(a_n\vee b_n \vee c_n)\notag\\
	&\rightarrow& \exp(\mathcal{N}(0,c^2)). \label{al_records_Burr_IV}
\end{eqnarray}

\Bin \textbf{Burr V}.\\

\Ni (i) The quantile function is expanded as follows.

 \begin{equation}
\frac{\pi }{2}-F^{-1}(1-u)=\left(\log \left( \frac{kr}{u}\right)\right)^{-1} - \frac{1}{2} \left(\log\left(\frac{kr}{u}\right)\right)^{-3}+O\left( \left(\log(1/u)\right)^{-5}\right),
\label{ResquantDefVa}
\end{equation}

\Bin i.e.,

\begin{equation}
\frac{\pi }{2}-F^{-1}(1-u)=\left(\log \left( \frac{kr}{u}\right)\right)^{-1}\biggr(1 - \frac{1}{2} \left(\log\left(\frac{kr}{u}\right)\right)^{-2}+O\left( \left(\log(1/u)\right)^{-4}\right)\biggr).
\label{ResquantDefVb}
\end{equation}

\Bin (ii) $F \in D(G_{\gamma})$, $\gamma=0$, $uep(F)=\pi/2$.\\

\Ni (iii) The asymptotic law of the record values $X^{(n)}$, $n\geq 1$, is given as follows.

\begin{eqnarray}
\sqrt{n}\biggr(\log\biggr(\frac{\pi}{2} - X^{(n)}\biggr)+\log n\biggr)&=&- S_{n}^{\ast} + O_{\mathbb{P}}(n^{-1/2})=- W_{n}^{\ast} + O_{\mathbb{P}}(c_n) \notag\\
&\rightarrow& \mathcal{N}(0,1).  \label{al_records_Burr_V}
\end{eqnarray}

\Bin We also have

\begin{eqnarray*}
\frac{(\log r + n)^{2} \biggr(X^{(n)}- \arctan \biggr(-\log\left\{\frac{(1-e^{-n})^{-1/r}-1}{k} \right\} \biggr)\biggr)}{\sqrt{n}}&=&S_{n}^{\ast} + O_{\mathbb{P}}(n^{-1/2})=W_n^\ast + O_{\mathbb{P}}(c_n)\notag\\
			&\rightarrow &\mathcal{N}(0,1). \ \ \  \text{\textit{(VAlt)}}
\end{eqnarray*}

\Bin \textbf{Burr VI}.\\

\Ni (i) The quantile function is expanded as follows.\\ 
 
\begin{equation}
F^{-1}(1-u)=\log 2+\log \log kr+\log \log \left( 1/u\right) +\frac{1}{4}%
\left( \log\log \frac{kr}{u}\right) ^{-2}+O\left( \log \log \left(
1/u\right)^{-3}\right), \ u<\frac{1}{e}.  \label{ResquantDefVI}
\end{equation}

\Ni (ii) $F \in D(G_{\gamma})$, $\gamma=0$, $ uep(F)=+\infty$.\\

\Ni (iii) The asymptotic law of the record values $X^{(n)}$, $n\geq 1$, is given as follows.

\begin{eqnarray}
\sqrt{n}\biggr(\frac{1}{(2n \log(kr))} \exp(X^{(n)})-1\biggr)&=& S_n^\ast + O_{\mathbb{P}}((\log n)^{-3})\notag\\
&=& W_n^\ast + O_{\mathbb{P}}((\log n)^{-3}) \notag\\
&\rightarrow & \mathcal{N}(0,1) \label{al_records_Burr_VI}
\end{eqnarray}

\Bin We also have

\begin{eqnarray*}
\frac{(\log r + n)^{2} \biggr(X^{(n)}-arcsinh \bigg( -\log \left(\frac{(1-e^{-n})^{-1/r}-1}{k} \right) \biggr)\biggr)}{\sqrt{n}}&=&S_{n}^{\ast} + O_{\mathbb{P}}(n^{-1/2}) \\
&=& W_n^\ast + O_{\mathbb{P}}(c_n)\notag\\
&\rightarrow &\mathcal{N}(0,1). \ \ \  \text{\textit{(VIAlt)}}
\end{eqnarray*}

\Bin \textbf{Burr VII}.\\

\Ni (i) The quantile function is expanded as follows.\\ 

\begin{equation}
F^{-1}\left( 1-u\right) =\log \sqrt{r}+\frac{1}{2}\log \left( 1/u\right) -%
\frac{1+r}{4r}u+O\left( u^{2}\right) .  \label{ResquantDefVII}
\end{equation}

\Bin (ii) Here $Z=\exp(X) \in D(G_{\gamma})$, $\gamma=1/2$, i.e., $X=\log Z$ and $uep(F)=+\infty $.\\
 
\Ni (iii) The asymptotic law of the record values $X^{(n)}$, $n\geq 1$, is given as follows.

\begin{eqnarray}
\frac{X^{(n)}-n/2 - \log \sqrt{r}}{\sqrt{n}}&=&(1/2) S_{n}^{\ast}+O_{\mathbb{P}}(d_n(\eta))= (1/2) W_{n}^{\ast}+ O_{\mathbb{P}}(c_n)\notag\\
																											&\rightarrow& \mathcal{N}(0, 1/4). \label{al_records_Burr_VII}
\end{eqnarray}

\Bin \textbf{Burr VIII}.\\

\Ni (i) The quantile function is expanded as follows.\\ 

\begin{equation}
F^{-1}(1-u) =\log (2r/\pi )+\log (1/u)-\frac{1-r}{2r}u+O(u^{2}). \label{ResquantDefVIII}
\end{equation}
 
\Bin (ii) Here $Z=\exp(X) \in D(G_{\gamma})$, $\gamma=1$, i.e., $X=\log Z$ and $uep(F)=+\infty$.\\ 

\Ni (iii) The asymptotic law of the record values $X^{(n)}$, $n\geq 1$, is given as follows.

\begin{eqnarray}
\frac{X^{(n)} - n - \log(2r/ \pi)}{\sqrt{n}}&=&S_{n}^{\ast}+O_{\mathbb{P}}(d_n(\eta))= W_{n}^{\ast}+ O_{\mathbb{P}}(c_n) \notag \\
&\rightarrow& \mathcal{N}(0,1). \label{al_records_Burr_VIII}
\end{eqnarray}

\Bin \textbf{Burr IX}.\\

\Ni (i) The quantile function is expanded as follows.\\

\begin{eqnarray*}
F^{-1}(1-u) &=&\left\{ 
\begin{array}{lll}
\frac{1}{r}\log \left( \frac{2}{uk}\right) -\left( \frac{2-k}{2r}\right)
u+O(u^{2}) & \text{if} & 0<r\leq 1/2 \label{ResquantDefIX-1}\\ 
&  &  \notag \\ 
\frac{1}{r}\log \left( \frac{2}{uk}\right) -\left( \frac{2-k}{2r}\right)
u+O(u^{1/r}) & \text{if} & r>1/2 \label{ResquantDefIX-2} 
\end{array} 
\right. 
\notag \\~~~~~~~~ \\
\end{eqnarray*}

\Ni (ii) Here $Z=\exp(X) \in D(G_{\gamma})$, $\gamma=1/r$, i.e., $X=\log Z$ and $uep(F)=+\infty$.\\

\Ni (iii) The asymptotic law of the record values $X^{(n)}$, $n\geq 1$, is given as follows for both sub-cases.

\begin{eqnarray}
\frac{X^{(n)} - \frac{n}{r} - \frac{1}{r} \log\left(\frac{2}{k}\right)}{\sqrt{n}}&=&\frac{1}{r}S_{n}^{\ast}+O_{\mathbb{P}}(d_n(\eta))=\frac{1}{r}W_{n}^{\ast}+ O_{\mathbb{P}}(c_n) \notag \\
	&\rightarrow& \mathcal{N}(0,r^{-2}). \label{al_records_Burr_IX}
\end{eqnarray}

\Ni \textbf{Burr X}.\\

\Ni (i) The quantile function is expanded as follows.\\ 
 
\begin{equation*}
F^{-1}(1-u)=(\log (1/u))^{1/2}\biggr\{1-\frac{r+1}{4r}\frac{u}{\log (1/u)}%
+O\left( \frac{u^{2}}{\log (1/u)}\right) \biggr\}. \label{ResquantDefX}
\end{equation*}

\Bin (ii) Here $X \in D(G_{\gamma})$, $\gamma=0$, $uep(F)=+\infty$.\\

\Ni (iii) The asymptotic law of the record values $X^{(n)}$, $n\geq 1$, is given as follows.

\begin{eqnarray}
\sqrt{n}\log \biggr(\frac{X^{(n)}}{\sqrt{n}}\biggr)&=& \frac{1}{2} S_{n}^{\ast} + O_{\mathbb{P}}(n^{-1/2})\notag\\
&=&\frac{1}{2} W_{n}^{\ast} + O_{\mathbb{P}}(c_n) \notag\\
&\rightarrow& \mathcal{N}(0,1/4). \label{al_records_Burr_X} 
\end{eqnarray}

\Bin We also have

\begin{eqnarray*}
\frac{(\log r + n)^{2} \biggr(X^{(n)}-\biggr(-\log\left((1-e^{-n})^{-1/r}-1\right)\biggr)^{1/2}\biggr)}{\sqrt{n}}&=&S_{n}^{\ast} + O_{\mathbb{P}}(n^{-1/2})=W_n^\ast + O_{\mathbb{P}}(c_n)\notag\\
	&\rightarrow &\mathcal{N}(0,1). \ \ \  \text{\textit{(XAlt)}}
\end{eqnarray*}

\Ni \textbf{Burr XI}.\\

\Ni (i) The quantile function is expanded as follows.\\ 
 
\begin{equation}
1 -F^{-1}\left( 1-u\right) =\alpha ^{-1/3}u^{1/3}\left( 1-%
\frac{\beta }{3\alpha }\alpha ^{-1/3}u^{2/3}+O\left( u^{4/9}\right) \right),
\label{ResquantDefXI}
\end{equation}

\Bin \Bin where $\alpha =\frac{(2\pi )^{2}}{6r},$ $\beta =-\frac{(2\pi )^{4}}{120r}$.\\

\Ni (ii) $F \in D(G_{\gamma})$, $\gamma=-1/3$, $uep(F)=1$.\\

\Ni (iii) The asymptotic law of the record values $X^{(n)}$, $n\geq 1$, is given as follows.

\begin{eqnarray}
\biggr(\alpha^{1/3} e^{n/3}\left(1 - X^{(n)}\right)\biggr)^{n^{-1/2}}&=& \exp \biggr(-\frac{1}{3} S_{n}^{\ast}\biggr)+O_{\mathbb{P}}(a_n \vee b_n ) \notag\\
&=& \exp \biggr(-\frac{1}{3} W_{n}^{\ast}\biggr)+ O_{\mathbb{P}}(a_n \vee b_n \vee c_n)\notag\\
																		&\rightarrow& \exp \biggr(\mathcal{N}(0,1/9) \biggr). \label{al_records_Burr_XI}
\end{eqnarray}

\Ni \textbf{Burr XII}.\\

\Ni (i) The quantile function is expanded as follows.\\ 
 
\begin{equation}
F^{-1}(1-u)=u^{-1/(rc)}\biggr(1-\frac{1}{c}u^{1/r}+\frac{1-c}{2c^{2}}%
u^{2/r}+O(u^{3/r})\biggr). \label{ResquantDefXII}
\end{equation}

\Bin (ii) $F \in D(G_{\gamma})$, $\gamma=1/(rc)$, $uep(F)=+\infty$.\\

\Ni (iii) The asymptotic law of the record values $X^{(n)}$, $n\geq 1$, is given as follows.

\begin{eqnarray}
\biggr(\frac{X^{(n)}}{e^{n/rc}}\biggr)^{n^{-1/2}}&=& \exp \biggr(\frac{1}{rc} S_{n}^{\ast} \biggr)+O_{\mathbb{P}}(a_n \vee b_n) \notag\\
&=& \exp \biggr(\frac{1}{rc} W_{n}^{\ast} \biggr) + O_{\mathbb{P}}(a_n \vee b_n \vee c_n)\notag\\
&\rightarrow& \exp \biggr(\mathcal{N}(0, (rc)^{-2}) \biggr). \label{al_records_Burr_XII}
\end{eqnarray}

\Bin \textbf{Distribution  (Xa)}.\\

\Ni (i) The quantile function is expanded as follows.\\ 

\begin{equation}
F^{-1}(1-u)= (\log(1/u))^{1/2} \biggr(1 + \frac{1-r}{4r} \frac{u}{\log(1/u)}
+ O\left(\frac{u^2}{\log(1/u)}\right)\biggr).  \label{ResquantDefXa} 
\end{equation}

\Bin (ii) Here $\in D(G_{\gamma}$, $\gamma=0$, $uep(F)=+\infty$.\\

\Ni (iii) The asymptotic law of the record values $X^{(n)}$, $n\geq 1$, is given as follows.

\begin{eqnarray}
\sqrt{n}\log \biggr(\frac{X^{(n)}}{\sqrt{n}}\biggr)&=& \frac{1}{2} S_{n}^{\ast} + O_{\mathbb{P}}(n^{-1/2})=\frac{1}{2} W_{n}^{\ast} + O_{\mathbb{P}}(c_n) \notag\\
&\rightarrow& \mathcal{N}(0,1/4). \label{al_records_Burr_Xa} 
\end{eqnarray}

\Bin We also have

\begin{eqnarray*}
\frac{(\log r + n)^{2} \biggr(X^{(n)}-\biggr(-\log\left((1-e^{-n})^{-1/r}-1\right)\biggr)^{1/2}\biggr)}{\sqrt{n}}&=&S_{n}{^\ast} + O_{\mathbb{P}}(n^{-1/2}) \notag\\
&=& W_n^\ast + O_{\mathbb{P}}(c_n)\notag\\
	&\rightarrow &\mathcal{N}(0,1). \ \ \  \text{\textit{(XaAlt)}}
\end{eqnarray*}

\end{theorem}

\Bin \textbf{Proofs}. As announced, we are going to provide the full proof of one case in which one the four arguments is applied. But for cases corresponding to 
Burr \textit{cdf}'s attracted to $D(G_0)$ with  $s(u)\rightarrow 0$ as $u\rightarrow 0$ (Burr V, Burr VI, Burr X, Burr Xa) , it is easier to draw the asymptotic law of record values directly from the quantile function expansions which are Formulas (\ref{al_records_Burr_V}), (\ref{al_records_Burr_VI}), (\ref{al_records_Burr_X}) and (\ref{al_records_Burr_Xa}) respectively. We begin by giving the outlines of the proofs using one of the three points (a,c,d) of Theorems \ref{theoALL1} and/or \ref{theoALL3}.\\

\Bin Next, we give the direct proofs for Burr \textit{cdf}'s attracted to $D(G_0)$ with  $s(u)\rightarrow 0$ as $u\rightarrow 0$.\\

\Bin However, in Appendix (A1), page \pageref{pageAppendixA1}, we will give alternative forms of the asymptotic laws of record values derived form Theorem \ref{extensionBelow} corresponding to Formulas \text{\textit{(VAlt)}}, \text{\textit{(VIAlt)}}, \text{\textit{(XAlt)}} and \text{\textit{(XaAlt)}} respectively.\\

\textbf{A - Proofs based on Direct applications of Theorems \ref{theoALL1} and/or \ref{theoALL3} and/or \ref{extensionBelow}}.\\

\Ni \textbf{Burr I}. We have $uep(F)=1$. Next, by Part (b) of Proposition \ref{criteria}, we have $F \in D(G_{\gamma})$ for $\gamma=-1$. The asymptotic law of the record values and the rates of convergence follow from Parts (c) of Theorems \ref{theoALL1} and \ref{theoALL3}.\\

\Ni \textbf{Burr II}. We have $uep(F)=+\infty$ and

$$
exp (F^{-1}(1-u))= r u^{-1} (1+o(u)), \ u \in ]0,1[.
$$

\Bin By Part (a) of Proposition \ref{criteria}, $\exp(X)$ with $ cdf$ \  $ F \in D(G_{\gamma})$ for $\gamma=1$. So, by Proposition \ref{proplo86}, $F \in D(G_{0})$.  The asymptotic law of the record values and the rates of convergence follow from Theorems \ref{theoALL1} and \ref{theoALL3}.\\

\Ni \textbf{Burr III}. We have $uep(F)=+\infty$. Next, by Part (a) of Proposition \ref{criteria}, we have $F \in D(G_{\gamma})$ for $\gamma=1/k$. The asymptotic law of the record values and the rates of convergence follow from Parts (a) of Theorems \ref{theoALL1} and \ref{theoALL3}.\\

\Ni \textbf{Burr IV}. We have $uep(F)=c$. Next, by Part (b) of Proposition \ref{criteria}, we have $F \in D(G_{\gamma})$ for $\gamma=-c$. The asymptotic law of the record values and the rates of convergence follow from Parts (b) of Theorems \ref{theoALL1} and \ref{theoALL3}.\\

\Ni \textbf{Burr V}. See Part B below .\\

\Ni \textbf{Burr VI}. See Part B below .\\

\Ni \textbf{Burr VII}. By Part (a) of Proposition \ref{criteria}, $\exp(X) $ with $ cdf$ \ $ F \in D(G_{\gamma})$ for $\gamma=1/2$. So, by Proposition \ref{proplo86}, $F \in D(G_{0})$.  The asymptotic law of the record values and the rates of convergence follow from Parts (b) in Theorems \ref{theoALL1} and \ref{theoALL3}.\\

\Ni \textbf{Burr VIII}. By Part (a) of Proposition \ref{criteria}, $\exp(X) $ with $ cdf$ \ $ F \in D(G_{\gamma})$ for $\gamma=1$. So, by Proposition \ref{proplo86}, $F \in D(G_{0})$.  The asymptotic law of the record values and the rates of convergence follow from Parts (b) in Theorems \ref{theoALL1} and \ref{theoALL3}.\\

\Ni \textbf{Burr IX}. By Part (a) of Proposition \ref{criteria}, $\exp(X) $ with $ cdf$ \ $ F \in D(G_{\gamma})$ for $\gamma=1/r$. So, by Proposition \ref{proplo86}, $F \in D(G_{0})$.  The asymptotic law of the record values and the rates of convergence follow from Parts (b) in Theorems \ref{theoALL1} and \ref{theoALL3}.\\

\Ni \textbf{Burr X}. See Part B below .\\

\Ni \textbf{Burr XI}. We have $uep(F)=1$. Next, by Part (b) of Proposition \ref{criteria}, we have $F \in D(G_{\gamma})$ for $\gamma=-1/3$. The asymptotic law of the record values and the rates of convergence follow from Parts (c) of Theorems \ref{theoALL1} and \ref{theoALL3}.\\

\Ni \textbf{Burr XII}. $uep(F)=+\infty$. Next, by Part (a) of Proposition \ref{criteria}, we have $F \in D(G_{\gamma})$ for $\gamma=1/(rc)$. The asymptotic law of the record values and the rates of convergence follow from Parts (a) of Theorems \ref{theoALL1} and \ref{theoALL3}.\\

\Ni \textbf{Burr Xa}. See Part B below.\\

\textbf{B - Proofs using direct methods}.\\

\Ni \textbf{Burr V}. We have to make a little effort to see that $F \in D(G_{0})$. We recall the quantile function (\ref{ResquantDefVa})

\begin{equation}
F^{-1}(1-u)=\frac{\pi }{2} -\left(\log \left( \frac{kr}{u}\right)\right)^{-1} + O\left( \left(\log(1/u)\right)^{-3}\right).
\end{equation}

\Bin Let $\ell(u)=(\log (kr/u)^{-1}$, $u\in ]0, \ u_0[$, $u_0 \in ]0, \ 1[$. We have $\ell^{\prime}(u)=(\log (kr/u))^{-2}/u$. So, for $s(u)=-\log (kr/u)^{-2}$, there existe a real number $c_0$ such that 

$$
\ell(u)=c_0 - \int_{u_0}^{u} \frac{s(t)}{t} \ dt.
$$

\Bin So $\ell(\circ)$ is slowly varying at $0$ and by the properties of slowly varying functions (See \cite{resnick},\cite{dehaan} and \cite{lo2018}), we have that for any $\lambda>0$, for all $u \in ]0, \ u_0[$,

$$
\lim_{u \rightarrow 0} \frac{\ell(\lambda u)-\ell(u)}{s(u)}=-\log \lambda.
$$

\Bin By remarking that for $g(u)=O\left( \left(\log(1/u)\right)^{-3}\right)$, we have $g(u)=O(s(u))$ as $u\rightarrow 0$. Hence, we may replace 
$\ell(u)$ by $F^{-1}(1-u)$ in the last limit to get

$$
\lim_{u \rightarrow 0} \frac{F^{-1}(1-\lambda u)-F^{-1}(1-u)}{s(u)}=- \log \lambda.
$$

\Bin So $F \in D(G_0)$.\\

\Ni To find the asymptotic law of record values, we can use Theorem \ref{extensionBelow} as in Appendix (A1), page \pageref{pageAppendixA1}. However, it is easier to use a direct method as below. From Expansion \ref{ResquantDefVb}, we have

\begin{eqnarray}
\frac{\pi}{2} -F^{-1}(1-u)=(\log(kr/u))^{-1} (1+O((\log(1/u))^{-2}), \ u \in ]0, \ 1[,
\end{eqnarray}

\Bin i.e.,

\begin{eqnarray*}
\log\biggr(\frac{\pi}{2} -F^{-1}(1-u)\biggr)=- (\log\log(kr/u)) + O((\log(1/u))^{-2}), \ u \in ]0, \ 1[.
\end{eqnarray*}

\Bin We have for $n\geq 1$,

\begin{eqnarray*}
\log\biggr(\frac{\pi}{2} - X^{(n)}\biggr)&=&\log \biggr(\frac{\pi}{2} - F^{-1}(1-e^{-S_{(n)}})\biggr)\\
&=&- \log \left(\log(kr) + S_{(n)}\right) + O_{\mathbb{P}}(n^{-2})\\
&=& - \log \biggr( S_{(n)} \left(1 +\frac{\log kr}{S_{(n)}}\right)\biggr) + O_{\mathbb{P}}(n^{-2})\\
&=& - \log S_{(n)} + O_{\mathbb{P}}(n^{-1}).\\
\end{eqnarray*}

\Bin Hence 

\begin{eqnarray*}
\log\biggr(\frac{\pi}{2} - X^{(n)}\biggr)+\log n&=& - \log (S_{(n)}/n) + O_{\mathbb{P}}(n^{-1})\\
&=&- \biggr(\log \biggr(1 + \frac{S_{(n)}}{n}-1\biggr)\biggr) + O_{\mathbb{P}}(n^{-1})\\
&=&- \frac{S_{(n)}-n}{n} + O_{\mathbb{P}}(n^{-1/2})+ O_{\mathbb{P}}(n^{-1})\\
&=&- \frac{S_{(n)}^{\ast}}{\sqrt{n}} + O_{\mathbb{P}}(n^{-1/2}). \ \square
\end{eqnarray*}

\Ni \textbf{Burr VI}. We proceed exactly as in the Proof related to Burr V with $s(u)=\biggr( \log(1/u)\biggr)^{-1}$, $0<u\leq u_{0} < 1 $, $c_0=\log\log(1/u_0)$,

$$
\log \biggr(\log(1/u)\biggr) =c_0 + \int_{u}^{u_0} \frac{s(t)}{t} \ dt.
$$

\Bin We get that $ F \in D(G_0)$. Let us use a direct method to find the asymptotic law of the record values. We apply the quantile function $V_{(n)}$ to get

\begin{eqnarray*}
\exp(X^{(n)})&=& 2\log(kr) S_{(n)} (1+O_{\mathbb{P}}((\log n)^{-3}))\\
\end{eqnarray*}

\Bin which leads to 

\begin{eqnarray*}
\frac{1}{2n\log(kr)} \exp(X^{(n)})&=& \frac{S_{(n)}}{n} (1+ O_{\mathbb{P}}((\log n)^{-3}))\\
&=& 1+ \biggr(\frac{S_{(n)}}{n}-1\biggr) (1+ O_{\mathbb{P}}((\log n)^{-3}))\\
&=& 1+\biggr( \frac{S_{(n)}- n}{n} \biggr)(1+ O_{\mathbb{P}}((\log n)^{-3})),
\end{eqnarray*}

\Bin and next

\begin{eqnarray*}
\sqrt{n} \biggr(\frac{1}{(2n)\log(kr)} \exp(X^{(n)})-1\biggr)&=& \frac{S_{(n)}-n}{\sqrt{n}}(1+ O_{\mathbb{P}}((\log n)^{-3}))\\
&=& S_{n}^\ast (1+ O_{\mathbb{P}}((\log n)^{-3})).
\end{eqnarray*}

\Bin Finally, we have

\begin{eqnarray*}
\sqrt{n} \biggr(\frac{1}{(2n)\log(kr)} \exp(X^{(n)})-1\biggr)&=& S_n^\ast + O_{\mathbb{P}}((\log n)^{-3}))\\
&=& W_n^\ast + O_{\mathbb{P}}((\log n)^{-3})) \\
& \rightarrow & \mathcal{N}(0,1) .\\
\end{eqnarray*}

\Ni \textbf{Burr X}. We have

\begin{equation*}
\log F^{-1}(1-u)=\frac{1}{2} \log \biggr(\log (1/u)\biggr)^{1/2}-\frac{r+1}{4r}\frac{u}{\log (1/u)}+O\left( \frac{u^{2}}{(\log (1/u))^{2}}\right). \label{ResquantDefX_02}
\end{equation*}

\Bin When applied to $V_{(n)}$, we get

$$
\log X^{(n)}= \frac{1}{2} \log S_{(n)} + O_{\mathbb{P}}((\log n)^{-2} d_n(\eta))
$$

\Bin and hence, by routine computations,

\begin{eqnarray*}
\log X^{(n)}- \frac{1}{2} \log n&=& \biggr(\frac{1}{2} \log S_{(n)} + O_{\mathbb{P}}((\log n)^{-2} d_n(\eta))\biggr)-\frac{1}{2}\log n \\
&=& \frac{1}{2} \log \frac{S_{(n)}}{n} + O_{\mathbb{P}}((\log n)^{-2} d_n(\eta))\\
&=& \frac{1}{2} \biggr(\left(\frac{S_{(n)}}{n}-1\right) + O_{\mathbb{P}}\left(\frac{S_{(n)}}{n}-1\right)\biggr)+O_{\mathbb{P}}((\log n)^{-2} d_n(\eta))\\
&=& \frac{1}{2} \frac{1}{\sqrt{n}} S_{n}^{\ast} + O_{\mathbb{P}}(n^{-1}).
\end{eqnarray*}

\Bin We conclude that

\begin{eqnarray*}
\sqrt{n}\log \biggr( \frac{X^{(n)}}{\sqrt{n}}\biggr)&=& \frac{1}{2} S_{n}^{\ast} + O_{\mathbb{P}}(n^{-1/2})=\frac{1}{2} W_{n}^{\ast} + O_{\mathbb{P}}(c_n)\\
&\rightarrow& \mathcal{N}(0,1/4).
\end{eqnarray*}

\Ni \textbf{Burr Xa}. We treat that case exactly as the case Burr X with the same final result,

\begin{equation}
F^{-1}(1-u)= (\log(1/u))^{1/2}+ \frac{1-r}{4r} \frac{u}{(\log(1/u))^{1/2}}
+ O\left(\frac{u^2}{(\log(1/u))^{1/2}}\right).  \label{ResquantDefXa_02} 
\end{equation}

\begin{eqnarray*}
\sqrt{n}\log \biggr(\frac{X^{(n)}}{\sqrt{n}}\biggr)&=& \frac{1}{2} S_{n}^{\ast} + O_{\mathbb{P}}(n^{-1/2})=\frac{1}{2} W_{n}^{\ast} + O_{\mathbb{P}}(c_n)\\
&\rightarrow& \mathcal{N}(0,1/4).
\end{eqnarray*}

$\blacksquare$\\

\Bin \textbf{Remark}.  By the way, once we have a second order extension of the quantile function, we may directly find the asymptotic law of the record value without applying arguments in Theorems \ref{theoALL1} and \ref{theoALL3}. However, those arguments, when put together, offer a unified approach to determine such asymptotic laws.\\

\section{Simulations} \label{sec_03}

Here we proceed to simulation studies of the asymptotic laws obtained. Since, we only want to illustrate how the results are for medium sizes, we will restric ourselves to two or three case in each extremal domain. The first issue we have to deal with concerns the sample sizes. In general, the sample size is fixed and the statistics using the generated sample are computed alongside related parameters. The situation is not the same with records theory.\\

\Ni Indeed, for a sample $X_1,\cdots,X_n$ of size $n\geq 1$, we study the $nr$-records. But, it is possible the sample does not have $nr$ records up to $n$ obervations. From \cite{ahsan-nev}, we have the following results. Let $N(n)$ be the number records in the sumple. The law of $N(n)$ is the sample is free-distribution. We have:

\begin{equation}
\mathbb{E}(N(n)) =(\log n) (1+o(1)) \ and \ \mathbb{V}ar(N(n))=(\log n) (1+o(1)), \label{nr-01}
\end{equation}

\begin{equation}
\frac{N(n)-\log n}{\sqrt{\log n}}\rightarrow \mathcal{N}(0,1), \label{nr-02}
\end{equation}

\Ni and

\begin{equation}
\liminf_{n\rightarrow +\infty} \frac{N(n)-\log n}{\sqrt{2\log n \log\log\log n}}=-11 \ and \ \limsup_{n\rightarrow +\infty} \frac{N(n)-\log n}{\sqrt{2\log n \log\log\log n}}=1. \label{nr-03}
\end{equation}

\Bin So, for $n$ fixed, we are not sure to have a fixed number of $nr$ records values. For example, by using the gaussian approximation, we have the following probability $p(3)$ of having less that $nr$ records values in a sample of size $n$ in Table \ref{Tab1} (see page \pageref{Tab1}).\\

\Bin So, while we want to have a powerful test for small samples, we should ensure that we have enough data to use a meaningful number of records. From Table \ref{Tab1}, we recommend to use the results for $n$ al least equal to $n=50$.\\

\begin{figure}
	\begin{centering}
		\includegraphics[scale=0.40]{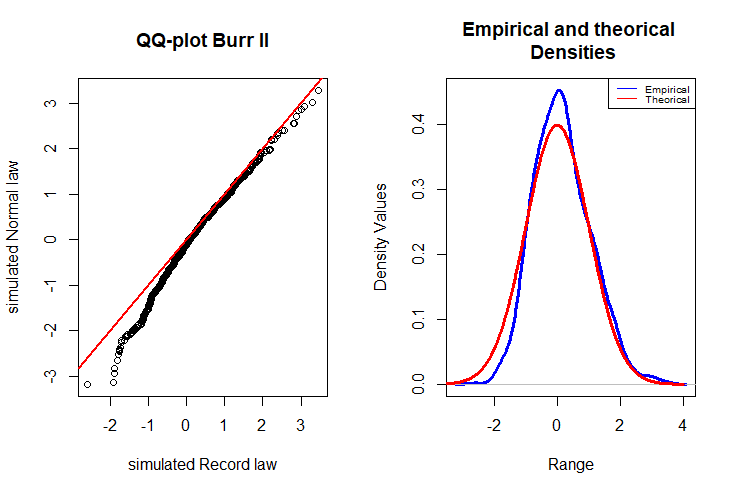}
		\par\end{centering}
	\begin{centering}
		\includegraphics[scale=0.40]{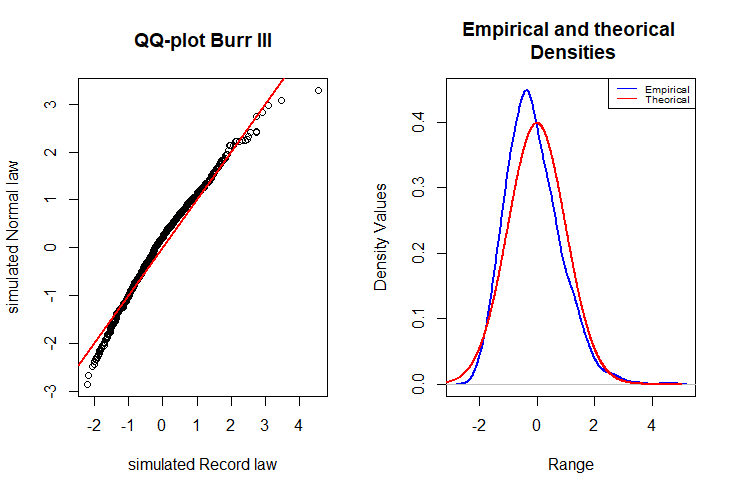}
		\par\end{centering}
	\caption{$\gamma>0:$ Burr II and Burr III}
	\label{fig1}
\end{figure}

\begin{figure}
	\begin{centering}
		\includegraphics[scale=0.4]{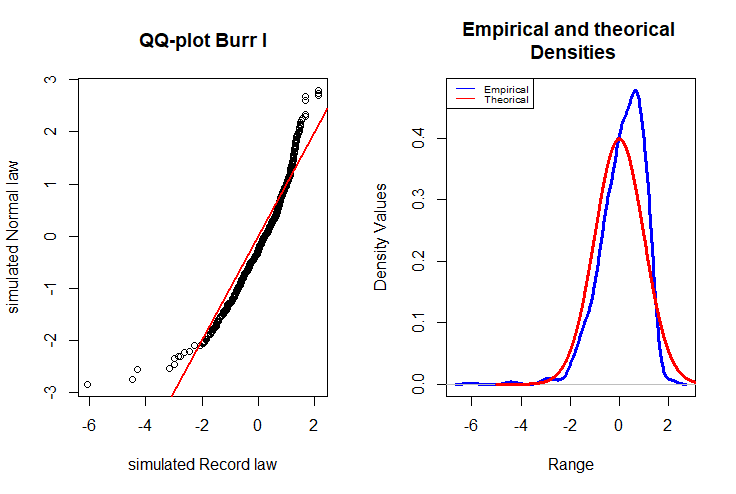}
		\par\end{centering}
	\begin{centering}
		\includegraphics[scale=0.4]{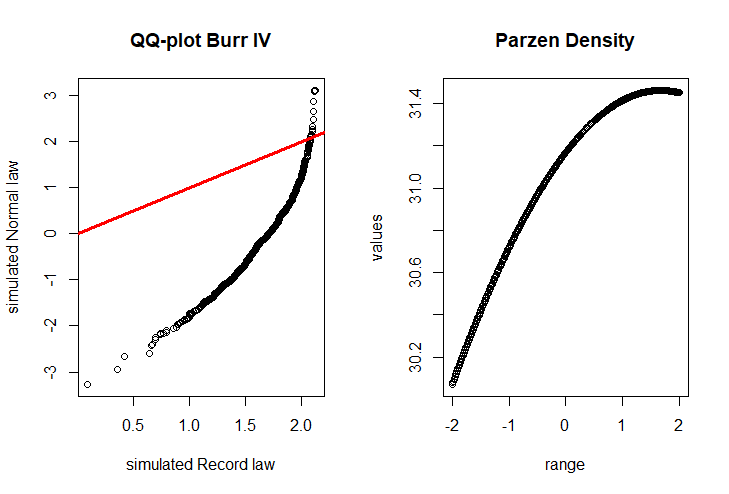}
		\par\end{centering}
	\caption{$\gamma<0:$ Burr I and Burr IV}
	\label{fig2}
	\end{figure}

\begin{figure}
	\begin{centering}
		\includegraphics[scale=0.4]{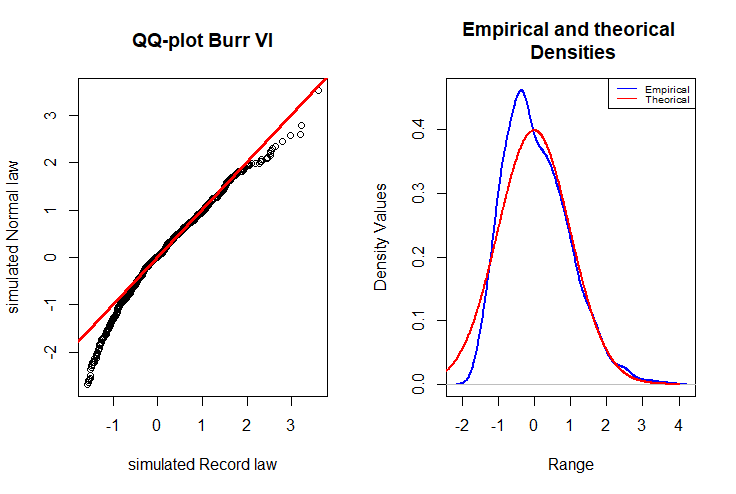}
		\par\end{centering}
	\begin{centering}
		\includegraphics[scale=0.4]{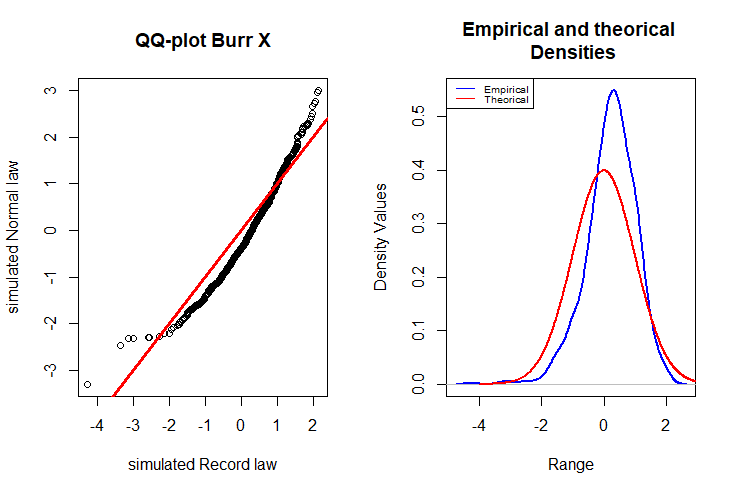}
		\par\end{centering}
	\caption{$\gamma=0:$ Burr VI and Burr X}
	\label{fig3}
\end{figure}

\Ni Now we are going to simulate the results on two for the cdf's of each extremal types: Burr II and Burr III for $\gamma>0$, Burr I and Burr IV for $\gamma<0$, Burr VI and XX for $\gamma=0$. For each of them we will compute the p-values of the asymptotic normality tests, and we display the qq-plots and the Parzen estimators of the pdf's of the centered and normalized records values. In Figures \ref{fig1} (for two Burr distributions in Type I), \ref{fig2} (for two Burr distributions in Type 
III), \ref{fig3} (for two Burr distributions in Type II), the \textit{qq-plots} and the Parzen graphs support our findings. Table \ref{Tab2} 
(in page \pageref{Tab2}) provides the p-values that validate our statistical tests.\\

\Ni \textbf{Description of the simulation works}. 

For each case, we proceed as follow
\begin{description}
	\item [{Step}] \textbf{1:} 
	
	Generate a $N-$sample of standard uniform law
	
	list the records obtained
	
	if the number of record values is enough (compared to the number of
	records fixed (nr) at beginnig)
	\begin{itemize}
		\item Compute the statistic test associated to stantard normal law, $Z[i]$
		\item Compute the proportion of absolute value of $Z$ greater than $1.96$
		(0.975-quantile of standard normal law)
	\end{itemize}
	\item [{Step}] \textbf{2:}
	
	Repeat \textbf{step 1}, $B=1000$ times
	
	report $P_{0}$ le mean of proportion obtained in the tries of \textbf{Step
		1}
	\item [{Step}] \textbf{3: Decision}
	
	If the value $P_{0}$ is less than $5\%$, we accept the normality
\end{description}

\Bin The analysis on the tables

\section{Second order expansions of Burr's quantile functions} \label{sec_04}

\bigskip \noindent \textbf{Quantile of Burr II distribution or parameter $r>0$} \label{quantBurrII}.\newline

\noindent Its support is $\mathcal{V}=\mathbb{R}$ and its \textit{cdf} is

\begin{equation*}
1-u=\biggr(1+e^{-x}\biggr)^{-r}, \ x\in \mathcal{V}, \ u \in ]0,1[. \ \ (\text{\textit{FII}})
\end{equation*}

\bigskip \noindent First, we will repeatedly need the following expansions,
as $u\rightarrow 0$,

\begin{eqnarray}
&&(1-u)^{-1/r}= 1 + \frac{u}{r} + \frac{r+1}{2r^2} u^2 + O(u^3)
\label{m-u-m-1overr} \\
&&(1-u)^{1/r}= 1 - \frac{u}{r} + \frac{1-r}{2r^2} u^2 + O(u^3).
\label{m-u-p-1overr}
\end{eqnarray}

\bigskip \noindent By applying Expansion \eqref{m-u-m-1overr} on (\textit{FII%
}), we get

\begin{eqnarray*}
-x &=&\log\biggr(\frac{u}{r} + \frac{r+1}{2r^2} u^2 + O(u^3)\biggr) \\
&=&\log\biggr(\frac{u}{r} \biggr(1+ \frac{r+1}{2r} u + O(u^2)\biggr)\biggr)
\\
&=&-\log r + \log u + \log\biggr(1+ \frac{r+1}{2r} u + O(u^2)\biggr). \\
\end{eqnarray*}

\bigskip \noindent Now we develop the logarithm in $v=\frac{r+1}{2r} u +
O(u^2)=O(u) \rightarrow 0$ at the first order, we get

\begin{equation*}
-x=-\log r + \log u + \frac{r+1}{2r} u + O(u^2),
\end{equation*}

\bigskip \noindent and we conclude

\begin{equation}
F^{-1}(1-u)=\log r + \log (1/u) - \frac{r+1}{2r} u + O(u^2).
\label{quantDefII}
\end{equation}

\bigskip \noindent \textbf{Quantile of Burr III distribution of
parameters $k>0$ and $r>0$} \label{quantBurrIII}.\newline

\noindent Its support is $\mathcal{V}=\mathbb{R}_+$ and its \textit{cdf} is

\begin{equation*}
1-u=\biggr(1 + x^{-k}\biggr)^{-r}, \ x\in \mathcal{V}, \ u \in ]0,1[. \ \ (%
\text{\textit{FIII}})
\end{equation*}

\bigskip \noindent By applying Expansion \eqref{m-u-m-1overr} on (\textit{%
FIII}), we get

\begin{eqnarray*}
x &=&\biggr(\frac{u}{r} \biggr(1+ \frac{r+1}{2r} u + O(u^2)\biggr)\biggr)%
^{-1/k} \\
&=& r^{1/k} u^{-1/k} \biggr(1+ \frac{r+1}{2r} u + O(u^2)\biggr)^{-1/k} \\
&=& r^{1/k} u^{-1/k} \biggr(1 - \frac{r+1}{2kr} u + O(u^2)\biggr). \\
\end{eqnarray*}

\bigskip \noindent We conclude

\begin{equation}
F^{-1}(1-u)= r^{1/k} u^{-1/k} \biggr(1 - \frac{r+1}{2kr} u + O(u^2)\biggr).
\label{quantDefIII}
\end{equation}

\bigskip \noindent \textbf{Quantile of Burr IV distribution of
parameters $c>0$ and $r>0$.} \label{quantBurrIV}.\newline

\noindent Its domain is $\mathcal{V}=[0,c]$ and its \textit{cdf} is given by

\begin{equation*}
F(x)=\left[1+\left(\frac{c-x}{x}\right)^{1/c}\right]^{-r}, \ x \in ]0,c].
\end{equation*}

\bigskip \noindent We have for $F(x)=1-u$ with $u \in [0;1[$,

\begin{equation*}
1-u=\left[1+\left(\frac{c-x}{x}\right)^{1/c}\right]^{-r}
\end{equation*}

\begin{eqnarray*}
\left(1-u\right)^{-1/r} &=&1+\left(\frac{c-x}{x}\right)^{1/c} \\
1+\frac{1}{r}u+ \frac{r+1}{2r^{2}}u^{2} +O\left(u^{3}\right)&=& 1+\left(%
\frac{c-x}{x}\right)^{1/c} \\
\frac{1}{r}u+ \frac{r+1}{2r^{2}}u^{2} +O\left(u^{3}\right)&=& \left(\frac{c-x%
}{x}\right)^{1/c} .\\
\end{eqnarray*}

\bigskip \noindent Hence, we have 
\begin{eqnarray*}
\left(\frac{c-x}{x}\right)^{1/c} &=&\frac{1}{r}u+ \frac{r+1}{2r^{2}}u^{2}
+O\left(u^{3}\right) \\
\left(\frac{c}{x}-1\right)^{1/c} &=& \frac{1}{r}u \left(1+\frac{r+1}{2r}u
+O\left(u^{2}\right)\right).
\end{eqnarray*}

\bigskip \noindent The last equation leads to, 

\begin{eqnarray*}
\frac{c}{x}-1 &=&\left[ \frac{1}{r}u\left( 1+\frac{r+1}{2r}u+O\left(
u^{2}\right) \right) \right] ^{c} \\
&=&r^{-c}u^{c}\left( 1+\frac{r+1}{2r}u+O\left( u^{2}\right) \right) ^{c} \\
&=&r^{-c}u^{c}\left( 1+\frac{c(r+1)}{2r}u+O\left( u^{2}\right) \right).
\end{eqnarray*}

\bigskip \noindent Then, we have 
\begin{eqnarray}
\frac{c}{x} &=&1+ r^{-c}u^{c} \left(1+\frac{c(r+1)}{2r}u
+O\left(u^{2}\right)\right).  \label{burrIVeq1}
\end{eqnarray}

\bigskip \noindent The Equation (\ref{burrIVeq1}), leads to 

\begin{eqnarray}
\frac{x}{c}=\left[ 1+r^{-c}u^{c}\left( 1+\frac{c(r+1)}{2r}u+O\left(
u^{2}\right) \right) \right] ^{-1}.  \label{burrIVeq2}
\end{eqnarray}

\bigskip \noindent Since $c>0 $, we have $r^{-c}u^{c} \left(1+\frac{c(r+1)}{2r}u +O\left(u^{2}\right)\right) \rightarrow 0$ as $u\rightarrow 0$.

\bigskip \noindent Equation (\ref{burrIVeq2}), leads to
 
\begin{equation*}
\frac{x}{c}=1-r^{-c}u^{c}\left( 1+\frac{c(r+1)}{2r}u+O\left( u^{2}\right)
\right) .
\end{equation*}

\bigskip \noindent That leads to, 

\begin{equation*}
x=c-cr^{-c}u^{c}\left( 1+\frac{c(r+1)}{2r}u+O\left( u^{2}\right) \right) .
\end{equation*}

\bigskip \noindent Finally, we have 

\begin{equation}
uep(F)-F^{-1}(1-u)=cr^{-c}u^{c}\left( 1+\frac{c(r+1)}{2r}u+O\left(
u^{2}\right) \right). \label{quantDefIV}
\end{equation}

\bigskip \noindent \textbf{Quantile of Burr V distribution of
parameters $k>0$ and $r>0$} \label{quantBurrV}.\newline

\noindent Its support is $\mathcal{V}=[-\pi/2, \ \pi/2]$ and its \textit{cdf}
is

\begin{equation*}
1-u=\biggr(1+k e^{-\tan x}\biggr)^{-r}, \ x\in \mathcal{V}, \ u \in ]0,1[. \
\ (\text{\textit{FV}})
\end{equation*}

\bigskip \noindent By applying Expansion \eqref{m-u-m-1overr} on (\textit{FV}%
), we get $u=1-F\left( x\right)$,

\begin{equation*}
\left( 1-u\right) ^{-1/r}=1+ke^{-\tan \left( x\right)}.
\end{equation*}

\bigskip \noindent By (\ref{m-u-m-1overr}) , we get

\begin{equation*}
1+ke^{-\tan \left( x\right) }=1+\frac{u}{r}+\frac{r-1}{2r^{2}}%
u^{2}+O\left(u^{3}\right).
\end{equation*}

\bigskip \noindent And so, we get

\begin{equation}
e^{-\tan \left( x\right) }=\frac{1}{k}\left( \frac{u}{r}+\frac{\left(
r-1\right) }{2r^{2}}u^{2}+O\left( u^{3}\right) \right) .  \label{eqR1}
\end{equation}

\bigskip \noindent Let us set

\begin{equation*}
h\left( x\right) =e^{-\tan \left( x\right) }, \ \ x\in ]-\frac{\pi }{2},%
\frac{\pi }{2}[.
\end{equation*}

\bigskip \noindent So, we get

\begin{equation*}
x=h^{-1}\left( \frac{1}{k}\left( \frac{u}{r}+\frac{\left( r-1\right) }{2r^{2}%
}u^{2}+O\left( u^{3}\right) \right) \right) .
\end{equation*}

\bigskip \noindent Let us find $h^{-1}$. We begin by remarking that

\begin{equation*}
\forall x\in ]-\frac{\pi }{2},\frac{\pi }{2}[, \ \tan \left( x\right) =\frac{%
1}{\tan \left( \frac{\pi }{2}-x\right)}.
\end{equation*}
\bigskip \noindent We set

\begin{equation*}
X=\frac{\pi }{2}-x.
\end{equation*}

\bigskip \noindent and remark that $X\rightarrow 0+$ as $x\rightarrow \left(
\pi /2\right) ^{-}$.\newline

\noindent We expand $\tan \left( X\right) $ as follows

\begin{equation*}
\tan \left( X\right) =X+\frac{X^{3}}{3}+\frac{2}{15}X^{5}+O\left(
X^{7}\right) .
\end{equation*}

\bigskip \noindent Hence,

\begin{eqnarray*}
\tan \left( x\right) &=&\frac{1}{\tan \left( X\right) }=X^{-1}\left( 1+\frac{%
X^{2}}{3}+\frac{2}{15}X^{4}+O\left( X^{6}\right) \right) ^{-1} \\
&=&X^{-1}\left( 1-\frac{X^{2}}{3}+\frac{7}{90}X^{4}+O\left( X^{6}\right)
\right).
\end{eqnarray*}

\bigskip \noindent We had already set

\begin{equation}
\tan (x)=Y=-\log y,as \  \ y \downarrow 0.  \label{eqR2}
\end{equation}

\bigskip \noindent So, we have

\begin{equation*}
Y=X^{-1}\left( 1-\frac{X^{2}}{3}+\frac{7}{90}X^{4}+O\left(
X^{6}\right)\right) .
\end{equation*}

\bigskip \noindent By the routine methods developed earlier, we have

\begin{equation*}
X=Y^{-1}\left( 1-\frac{1}{3} Y^{-2}+O\left( Y^{-4}\right) \right).
\end{equation*}

\bigskip \noindent So

\begin{equation}
\frac{\pi }{2}-x=\left(\log \left( 1/y\right)\right) ^{-1}\left( 1-\frac{\left(\log\left(1/y\right)\right) ^{-2}}{3}+O\left( \log \left( 1/y\right) ^{-4}\right)
\right).  \label{eqR3}
\end{equation}

\bigskip \noindent By formula \eqref{eqR2} , we have

\begin{equation*}
\tan x=-\log y\Longleftrightarrow e^{-\tan x}=y\Longleftrightarrow
h\left(x\right) =y.
\end{equation*}

\bigskip \noindent But from Formulas \eqref{eqR1} and \eqref{eqR2}, we may
take

\begin{equation*}
y=:y\left( u\right) =\frac{u}{kr}\left( 1+\frac{r+1}{2r}u+O\left(u^{2}%
\right) \right) .
\end{equation*}

\bigskip \noindent Hence, Formula \eqref{eqR3} becomes

\begin{eqnarray*}
\frac{\pi }{2}-x &=&\log \left( 1/y\left( u\right) \right) ^{-1}\left( 1-%
\frac{\log \left( 1/y\left( u\right) \right) ^{-2}}{3}+O\left( \log \left(
1/y\left( u\right) \right) ^{-4}\right) \right) \\
&=&\log \left( 1/y\left( u\right) \right) ^{-1}-\frac{\log \left( 1/y\left(
u\right) \right) ^{-3}}{2}+O\left( \log \left( 1/y\left( u\right) \right)
\right) ^{-5}.
\end{eqnarray*}

\bigskip \noindent But

\begin{eqnarray*}
\log \left( \frac{1}{y\left( u\right) }\right) &=&-\log \left( y\left(
u\right) \right) \\
&=&\log \left( \frac{kr}{u}\right) -\log \left[ 1+\frac{r+1}{2r}u+O\left(
u^{2}\right) \right] ^{-1} \\
&=&\log \left( \frac{kr}{u}\right) \left( 1-\frac{r+1}{2r}\frac{u}{\log 
\frac{kr}{u}}+O\left( \frac{u^{2}}{\log \frac{kr}{u}}\right) \right) .
\end{eqnarray*}

\bigskip \noindent Then

\begin{eqnarray*}
\left( \log \frac{1}{y\left( u\right) }\right) ^{-1} &=&\left( \log \frac{kr%
}{u}\right) ^{-1}\left(1+\frac{r+1}{2r}\frac{u}{\log \frac{kr}{u}}+O\left( 
\frac{u^{2}}{\log \frac{kr}{u}}\right) \right) \\
&=&\left( \log \frac{kr}{u}\right) ^{-1}+\frac{r+1}{2r}\frac{u}{\left( \log 
\frac{kr}{u}\right) ^{2}}+O\left( \left( \frac{u}{\log \frac{kr}{u}}%
\right)^{2}\right).
\end{eqnarray*}

\bigskip \noindent Finally, we have

\begin{equation}
\frac{\pi }{2}-F^{-1}(1-u)=\left(\log \left( \frac{kr}{u}\right)\right)^{-1} - \frac{1}{2} \left(\log\left(\frac{kr}{u}\right)\right)^{-3}+O\left( \left(\log(1/u)\right)^{-5}\right).
\label{quantDefV}
\end{equation}



\bigskip \noindent \textbf{Quantile of Burr VI distribution of
parameters $k>0$ and $r>0$} \label{quantBurrVI}.\newline

\noindent Its support is $\mathcal{V}=\mathbb{R}$ and its \textit{cdf} is

\begin{equation*}
1-u=\biggr(1+ k e^{- \sinh x}\biggr)^{-r}, \ x\in \mathcal{V}, \ u \in
]0,1[. \ \ (FVI)
\end{equation*}

\bigskip \noindent By applying Expansion \eqref{m-u-m-1overr} to Formula
(FVI), we have

\begin{equation*}
e^{-\sinh \left( x\right) }=\frac{1}{k}\left( \frac{u}{r}+\frac{r+1}{2r}%
u^{2}+O\left( u^{3}\right)\right).
\end{equation*}

\bigskip \noindent Thus

\begin{equation*}
y=e^{-\sinh \left( x\right) }\Longleftrightarrow \sinh \left( x\right)=-\log
y=:Y.
\end{equation*}

\bigskip \noindent Then

\begin{eqnarray*}
\sinh \left( x\right) &=&\frac{e^{x}-e^{-x}}{2}=Y \\
&\Longleftrightarrow &e^{2x}-2Ye^{x}-1=0.
\end{eqnarray*}

\bigskip \noindent \bigskip \noindent This equation has two solutions:

\begin{equation*}
e^{x}=Y-\sqrt{Y^{2}+1}(i),\text{ or }e^{x}=Y+\sqrt{Y^{2}+1} \ \ (ii).
\end{equation*}

\bigskip \noindent The solution (i) is impossible since for $Y\geq 0,Y-\sqrt{%
Y^{2}+1}\leq 0$ and so, $e^{x}\neq Y-\sqrt{Y^{2}+1}$ for any $x\in \mathbb{R}
$. We keep Solution $(ii)$. Hence

\begin{eqnarray}
x &=&\log Y+\log \left( 1+\left( 1+y^{-2}\right) ^{\frac{1}{2}}\right) \notag \\
&=&\log Y+\log \left( 2+\frac{1}{2}Y^{-2}-\frac{1}{8}Y^{-4}+O\left(Y^{-6}%
\right) \right) \notag \\
&=&\log Y+\log 2+\frac{1}{4}Y^{-2}-\frac{3}{32}Y^{-4}+O\left( Y^{-6}\right).
\label{eqR4}
\end{eqnarray}

\bigskip \noindent By (\ref{eqR4}), we have

\begin{eqnarray*}
y &=&y\left( u\right) =\frac{1}{k}\left( \frac{u}{r}+\frac{r+1}{2r^{2}}%
u^{2}+O\left( u^{3}\right) \right) \\
&=&\frac{u}{kr}\left( 1+\frac{r+1}{2r}u+O\left( u^{2}\right) \right).
\end{eqnarray*}

\bigskip \noindent So

\begin{eqnarray*}
-\log y\left( u\right) &=&\log \frac{kr}{u}-\frac{r+1}{2r}%
u+O\left(u^{2}\right) \\
&=&\log \frac{kr}{u}\left( 1-\frac{r+1}{2r}\frac{u}{\log \frac{kr}{u}}%
+O\left( \frac{u^{2}}{\log \frac{kr}{u}}\right) \right) .
\end{eqnarray*}

\bigskip \noindent Hence

\begin{eqnarray}
\log Y=\log \log \frac{kr}{u}-\frac{r+1}{2r}\left( \frac{u}{\log \frac{kr}{u}%
}\right) +O\left( \frac{u}{\log 1/u}\right).  \label{eqR5}
\end{eqnarray}

\bigskip \noindent By (\ref{eqR5}), we have

\begin{equation*}
Y^{-\alpha }=\left( \log \frac{kr}{u}\right) ^{-\alpha }\left( 1+\frac{%
\alpha \left( r+1\right) }{2r}\frac{u}{\log\frac{kr}{u}}+O\left( \frac{u}{%
\log \frac{kr}{u}}\right) ^{2}\right) .
\end{equation*}

\bigskip \noindent Finally, by (\ref{eqR4}),

\begin{equation}
F^{-1}(1-u)=\log 2+\log \log kr+\log \log \left( 1/u\right) +\frac{1}{4}%
\left( \log\log \frac{kr}{u}\right) ^{-2}+O\left( \log \log \left(
1/u\right)^{-3}\right), \ u<\frac{1}{e}.  \label{quantDefVI}
\end{equation}

\bigskip \noindent \textbf{Quantile of Burr VII distribution of
parameter $r>0$} \label{quantBurrVII}.\newline

\noindent Its support is $\mathcal{V}=\mathbb{R}$ and its \textit{cdf} is

\begin{equation*}
1-u=2^{r}\biggr(1+ \tanh x \biggr)^{r}, \ x\in \mathcal{V}, \ u \in ]0,1[. \
\ (\text{\textit{FVII}})
\end{equation*}

\bigskip \noindent By applying Expansion \eqref{m-u-p-1overr} on (\text{%
\textit{FVII}}), we get

\begin{equation*}
\tanh \left( x\right) =2\left( 1-u\right) ^{1/r}-1=1-\frac{2}{r}u+\frac{1-r}{%
r^{2}}u^{2}+O\left( u^{3}\right) =:y .
\end{equation*}

\bigskip \noindent So, we have

\begin{equation*}
\frac{e^{x}-e^{-x}}{e^{x}+e^{x}}=\frac{e^{2x}-1}{e^{2x}+1}=y\in ]-1,1[.
\end{equation*}

\bigskip \noindent Then,

\begin{equation*}
e^{2x}=\frac{y+1}{y-1}\text{ ,\ \ }\ y\in ]-1,1[,\ x\in \mathbb{R}.
\end{equation*}

\bigskip \noindent In the formula above, $x\uparrow +\infty$ as $%
y\uparrow 1$. So

\begin{eqnarray*}
x &=&\frac{1}{2}\log \frac{y+1}{y-1} \\
&=&\frac{1}{2}\left[ \log \left( 1+y\right) -\log \left( 1-y\right) \right]
\\
&=&\frac{1}{2}\left[ \log \left( 2-\frac{2}{r}u+\frac{1-r}{r^{2}}%
u^{2}+O\left( u^{3}\right) \right) -\log \left( \frac{2}{r}u\left( 1-\frac{%
1-r}{2r}u+O\left( u^{2}\right) \right) \right) \right] \\
&=&\frac{1}{2}\left[ \left( \log 2-\frac{1}{r}u+\frac{1-r}{2r^{2}}u^{2}-%
\frac{1}{2}\frac{1}{r^{2}}u^{2}+O\left( u^{3}\right) \right) -\left( \log 
\frac{2}{r}u-\frac{1-r}{2r}u+O\left( u^{2}\right) \right) \right] .
\end{eqnarray*}

\bigskip \noindent So we have,

\begin{equation*}
2x=\log r+\log \left( 1/u\right) -\frac{1+r}{2r}u+O\left( u^{2}\right).
\end{equation*}

\bigskip \noindent Thus,

\begin{equation*}
x=\log \sqrt{r}+\frac{1}{2}\log \left( 1/u\right) -\frac{1+r}{4r}u+O\left(
u^{2}\right).
\end{equation*}

\bigskip \noindent Hence,

\begin{equation}
F^{-1}\left( 1-u\right) =\log \sqrt{r}+\frac{1}{2}\log \left( 1/u\right) -%
\frac{1+r}{4r}u+O\left( u^{2}\right) .  \label{quantDefVII}
\end{equation}

\bigskip \noindent \textbf{Quantile of Burr VIII distribution of parameter $%
r>0$} \label{quantBurrVIII}.\newline

\noindent Its support is $\mathcal{V}=\mathbb{R}$ and its \textit{cdf} is

\begin{equation*}
1-u=\biggr(\frac{2}{\pi} \arctan (e^{x})\biggr)^{r}, \ x\in \mathcal{V}, \ u
\in ]0,1[. \ \ (\text{\textit{FVIII}})
\end{equation*}

\bigskip \noindent \textbf{Case $r\neq 1$}.\newline

\bigskip \noindent By applying Expansion \eqref{m-u-p-1overr} on (\textit{%
FVIII}), we get

\begin{eqnarray*}
x &=& \log\biggr[\tan \biggr(\frac{\pi}{2}\biggr\{1 - \frac{u}{r} + \frac{1-r%
}{2r^2} u^2 + O(u^3)\biggr\} \biggr)\biggr]. \\
\end{eqnarray*}

\bigskip \noindent From there, we use the property that $\tan(\pi/2 -
u)=1/\tan(u) $ for $u$ positive and small, to have

\begin{eqnarray*}
x &=&-\log \biggr[\tan \biggr(\biggr\{\frac{u}{r}+\frac{1-r}{2r^{2}}%
u^{2}+O(u^{3})\biggr\}\biggr)\biggr]. \\
&&
\end{eqnarray*}

\bigskip \noindent Next, using expansion $\tan(v)=v + v ^{3}/3 + 2 v^{5}/15
+ O(v^7)$ as $v\rightarrow 0$ but restricting to the first order, we have

\begin{eqnarray*}
x =&-\log \biggr[\frac{\pi }{2r}u\biggr(1-\frac{1-r}{2r}u+O(u^{2})\biggr)\biggr]. 
\end{eqnarray*}

\bigskip \noindent Finally, we have

\begin{equation}
F^{-1}(1-u) =\log (2r/\pi )+\log (1/u)-\frac{1-r}{2r}u+O(u^{2}). \label{quantDefVIII}
\end{equation}

\bigskip \noindent \textbf{Quantile of Burr IX distribution of parameters $%
k>0 $ and $r>0$} \label{quantBurrIX}.\newline

\noindent Its support is $\mathcal{V}=\mathbb{R}$ and its \textit{cdf} is

\begin{eqnarray*}
1-u &=&1-\biggr(\frac{2}{2+k\left( \left( 1+e^{x}\right) ^{r}-1\right) }%
\biggr),\ x\in \mathcal{V},\ u\in ]0,1[.\ \ (\text{\textit{FIX}}) \\
&&
\end{eqnarray*}

\bigskip \noindent We get

\begin{eqnarray*}
\frac{2}{u} &=&2+k\left( \left( 1+e^{x}\right) ^{r}-1\right),
\end{eqnarray*}

\bigskip \noindent which leads to

\begin{eqnarray*}
\left( 1+e^{x}\right) ^{r}-1 &=&\frac{1}{k}\left( \frac{2}{u}-2\right),
\end{eqnarray*}

\bigskip \noindent that is

\begin{eqnarray*}
\left( 1+e^{x}\right) ^{r} &=&\frac{2u^{-1}}{k}(1-u)+1 \\
&=&\frac{2u^{-1}}{k}\left( 1-u+\frac{ku}{2}\right) \\
&=&\frac{2u^{-1}}{k}\left( 1-\frac{2-k}{2}u\right) .
\end{eqnarray*}

\bigskip \noindent So, we have by expanding

\begin{eqnarray*}
1+e^{x} &=&\left( \frac{2u^{-1}}{k}\right) ^{1/r}\left( 1-\frac{2-k}{2r}%
u+O(u^{2})\right).
\end{eqnarray*}%

\Bin Then we have
 
\begin{eqnarray*}
e^{x} &=&\left( \frac{2u^{-1}}{k}\right) ^{1/r}\left( 1-\left( \frac{2-k}{2r}%
\right) u-\left( \frac{ku}{2}\right) ^{1/r}+O(u^{2})\right) \\
\end{eqnarray*}

\Bin and 

\begin{eqnarray*}
x &=&\frac{1}{r}\log \left( \frac{2u^{-1}}{k}\right) -\left( \frac{2-k}{2r}%
\right) u-\left( \frac{ku}{2}\right) ^{1/r}+O(u^{2}).
\end{eqnarray*}

\Bin Now, we conclude, that the quatile depends on the value of $r>0$, as follows.

\begin{eqnarray}
F^{-1}(1-u) &=&\left\{ 
\begin{array}{lll}
\frac{1}{r}\log \left( \frac{2}{uk}\right) -\left( \frac{2-k}{2r}\right)
u+O(u^{2}) & \text{if} & 0<r\leq 1/2 \label{quantDefIX-1}\\ 
&  &  \notag \\ 
\frac{1}{r}\log \left( \frac{2}{uk}\right) -\left( \frac{2-k}{2r}\right)
u+O(u^{1/r}) & \text{if} & r>1/2 .\label{quantDefIX-2}
\end{array}
\right. 
\end{eqnarray}

\Bin 
\bigskip \noindent \textbf{Quantile of Burr X distribution of parameter $%
r>0$} \label{quantBurrX}.\newline

\noindent Its support is $\mathcal{V}=\mathbb{R}_+$ and its \textit{cdf} is

\begin{equation*}
1-u=\biggr(1+e^{-x^2}\biggr)^{-r}, \ x\in \mathcal{V}, \ u \in ]0,1[. \ \ (%
\text{\textit{FXa}})
\end{equation*}

\bigskip \noindent By applying Expansion \eqref{m-u-m-1overr} on (\textit{FXa%
}), we get

\begin{eqnarray*}
e^{-x^2} &=&\frac{u}{r} \biggr(1 + \frac{r+1}{2r} u + O(u^2))\biggr) \\
\end{eqnarray*}

\bigskip \noindent and

\begin{eqnarray*}
-x^2 &=& -\log r + \log u + \frac{r+1}{2r} u + O(u^2)) \\
&=&- \biggr(\log(1/u)\biggr\{1 - \frac{r+1}{2r} \frac{u}{\log(1/u)} + \frac{%
\log r}{\log(1/u)}+ O\left(\frac{u^2}{\log(1/u)}\right)\biggr\}\biggr).
\end{eqnarray*}

\bigskip \noindent So by expanding the latter line at the power $1/2$, we get

\begin{eqnarray*}
x &=&\biggr(\log (1/u)\biggr\{1-\frac{r+1}{2r}\frac{u}{\log (1/u)}+\frac{%
\log r}{\log (1/u)}+O\left( \frac{u^{2}}{\log (1/u)}\right) \biggr\}\biggr)%
^{1/2} \\
&=&(\log (1/u))^{1/2}\biggr\{1-\frac{r+1}{4r}\frac{u}{\log (1/u)}+O\left( 
\frac{u^{2}}{\log (1/u)}\right) \biggr\}. \\
&&
\end{eqnarray*}

\Bin Finally, we have

\begin{equation*}
F^{-1}(1-u)=(\log (1/u))^{1/2}\biggr\{1-\frac{r+1}{4r}\frac{u}{\log (1/u)}%
+O\left( \frac{u^{2}}{\log (1/u)}\right) \biggr\}. \label{quantDefX}
\end{equation*}

\bigskip \noindent \textbf{Quantile of Burr XI distribution of parameter $%
r>0 $} \label{quantBurrXI}.\newline

\noindent Its support is $\mathcal{V}=[0,1]$ and its \textit{cdf} is

\begin{equation*}
F(x)=\biggr(x - \frac{1}{2\pi} \sin(2\pi x)\biggr)^{r}, \ x\in \mathcal{V},
\ u \in ]0,1[. \ \ (\text{\textit{FXI}})
\end{equation*}

\bigskip \noindent Let us set $g(x)=\sin(2\pi x)$ of $x\in [0,1]$. The first
seven derivatives (at left of $x=1$) are

\begin{equation*}
g^{\prime }(x)=2\pi \cos (2\pi x),\ \ g^{\prime \prime }(x)=-(2\pi )^{2}\sin
(2\pi x)\ ,\ g^{(3)}(x)=-(2\pi )^{3}\cos (2\pi x),
\end{equation*}

\begin{equation*}
g^{(4)}(x)=+(2\pi )^{4}\sin (2\pi x),\ \ g^{(5)}(x)=(2\pi )^{5}\cos (2\pi x),
\end{equation*}

\begin{equation*}
g^{(6)}(x)=-(2\pi )^{6}\sin (2\pi x)\text{ \ and \ }g^{(7)}(x)=-(2\pi
)^{6}\cos (2\pi x).
\end{equation*}

\bigskip \noindent The even derivatives vanish at $x=1$ and the odd
derivatives take values $(-1)^{k}(2\pi)^{2k+1}$, $k\geq 0$. Thus, $ g $ is
expanded at $x=1$ as follows.

\begin{equation*}
\sin(2\pi x)= (2\pi) (x-1) - \frac{(2\pi)^{3}}{6} (x-1)^3 + \frac{(2\pi)^{5}%
}{5!} (x-1)^5 - \frac{(2\pi)^{7}}{7!} (x-1)^7 + O\left((x-1)^9)\right).
\end{equation*}

\noindent So, we have

\begin{eqnarray*}
F(x) &=&\biggr(x-(x-1)+\frac{(2\pi )^{2}}{6}(x-1)^{3}-\frac{(2\pi )^{4}}{5!}%
(x-1)^{5}+\frac{(2\pi )^{6}}{7!}(x-1)^{7}+O\left( (x-1)^{9})\right) \biggr)%
^{r} \\
&=&\biggr(1+\frac{(2\pi )^{2}}{6}(x-1)^{3}-\frac{(2\pi )^{4}}{5!}(x-1)^{5}+%
\frac{(2\pi )^{6}}{7!}(x-1)^{7}+O\left( (x-1)^{9})\right) \biggr)^{r}.
\end{eqnarray*}%

\Bin We set  $v=\frac{(2\pi )^{2}}{6}(x-1)^{3}-\frac{(2\pi )^{4}}{5!}(x-1)^{5}+%
\frac{(2\pi )^{6}}{7!}(x-1)^{7}+O\left( (x-1)^{9})\right) $ and the \textit{cdf}
becomes   

\begin{eqnarray*}
F(x) &=&(1+v)^{r} \\
&=&1+rv+\frac{r(r-1)}{2}v^{2}+O(v^{3}) \\
&=&1+\frac{(2\pi )^{2}}{6r}(x-1)^{3}-\frac{(2\pi )^{4}}{120r}(x-1)^{5}+\frac{%
r(r-1)}{2}\frac{(2\pi )^{4}}{36}(x-1)^{6}+O\left( (x-1)^{7})\right).
\end{eqnarray*}

\bigskip \noindent Hence

\begin{eqnarray*}
1-F(x) &=&\frac{(2\pi )^{2}}{6r}(1-x)^{3}-\frac{(2\pi )^{4}}{120r}(1-x)^{5}-%
\frac{r(r-1)}{2}\frac{(2\pi )^{4}}{36}(x-1)^{6}+O\left( (x-1)^{7})\right) \\
&=&\frac{(2\pi )^{2}}{6r}(1-x)^{3}-\frac{(2\pi )^{4}}{120r}(1-x)^{5}+O\left(
(x-1)^{6}\right) \\
&=&\alpha X^{3}+\beta X^{5}+O\left( X^{6}\right) .
\end{eqnarray*}

\Bin where $\alpha =\frac{(2\pi )^{2}}{6r},$ $\beta =-\frac{(2\pi )^{4}}{120r}$ and $X=1-x.$ We set $u=1-F(x)$ and we have the following 

\begin{eqnarray*}
u &=&\alpha X^{3}+\beta X^{5}+O\left( X^{6}\right) \\
&=&\alpha X^{3}\left( 1+\frac{\beta }{\alpha }X^{5}+O\left( X^{3}\right)
\right).
\end{eqnarray*}

\Bin By the same method used previously, 

\begin{equation*}
X=\alpha ^{-1/3}u^{1/3}\left( 1-\frac{\beta }{3\alpha }\alpha
^{-1/3}u^{2/3}+O\left( u^{4/9}\right) \right) .
\end{equation*}

\Bin That leads to, 

\begin{equation*}
1-x=\alpha ^{-1/3}u^{1/3}\left( 1-\frac{\beta }{3\alpha }\alpha
^{-1/3}u^{2/3}+O\left( u^{4/9}\right) \right) .
\end{equation*}

\Bin Hence

\begin{equation}
uep\left( F\right) -F^{-1}\left( 1-u\right) =\alpha ^{-1/3}u^{1/3}\left( 1-%
\frac{\beta }{3\alpha }\alpha ^{-1/3}u^{2/3}+O\left( u^{4/9}\right) \right) .
\label{quantDefXI}
\end{equation}

\bigskip \noindent \textbf{Quantile of Burr XII distribution of parameters $%
c>0$ and $r>0$} \label{quantBurrXII}.\newline

\noindent Its support is $\mathcal{V}=\mathbb{R}_+$ and its \textit{cdf} is

\begin{equation*}
1-u=1 - \biggr(1+ x^c\biggr)^{-r}, \ x\in \mathcal{V}, \ u \in ]0,1[. \ \ (%
\text{\textit{FXII}})
\end{equation*}

\bigskip \noindent We have

\begin{eqnarray*}
x &=&u^{-1/(rc)}\biggr(1-u^{1/r}\biggr)^{1/c} \\
&=&u^{-1/(rc)}\biggr(1-\frac{1}{c}u^{1/r}+\frac{1-c}{2c^{2}}%
u^{2/r}+O(u^{3/r})\biggr). \\
&&
\end{eqnarray*}

\Bin Finally, we have

\begin{equation}
F^{-1}(1-u)=u^{-1/(rc)}\biggr(1-\frac{1}{c}u^{1/r}+\frac{1-c}{2c^{2}}%
u^{2/r}+O(u^{3/r})\biggr). \label{quantDefXII}
\end{equation}


\bigskip \noindent \textbf{Quantile of distribution (Xa) of parameter $%
r>0$} \label{quantBurrXa}.\newline


\noindent Its support is $\mathcal{V}=\mathbb{R}_+$ and its \textit{cdf} is

\begin{equation*}
1-u=\biggr(1-e^{-x^2}\biggr)^{r}, \ x\in \mathcal{V}, \ u \in ]0,1[. \ \ (%
\text{\textit{Fx}})
\end{equation*}

\bigskip \noindent We have

\begin{equation*}
x^2= -\log\biggr( 1 - (1-u)^{-1/r}\biggr)
\end{equation*}

\bigskip \noindent By using the computations as defined in Burr II's case, we get

\begin{equation}
F^{-1}(1-u)= (\log(1/u))^{1/2} \biggr(1 + \frac{1-r}{4r} \frac{u}{\log(1/u)}
+ O\left(\frac{u^2}{\log(1/u)}\right)\biggr).  \label{quantDefXa} 
\end{equation}

\section{Conclusion} \label{sec_05}

\newpage

\Bin \textbf{Appendix (A1) : Applying Theorem \ref{extensionBelow} for Burr V, VI, X and Xa}. \label{pageAppendixA1}\\

\Ni We compute $s(u)=-u(F^{-1}(1-u))^{\prime}$ for each of these four cases and remark that $s(u)\rightarrow 0$ as $u\rightarrow 0$ and $ s(\circ) $ is slowly varying at zero. We consider the expression of $F(\circ)$ and $F^{-1}(\circ)$ for Burr V, VI, X and Xa in Table \ref{tab01} and find the following expressions of $s(\circ)$ :

\begin{eqnarray*}
(Burr \ V)&:& F^{-1}(1-u)= arctan \biggr(-\log\left\{\frac{(1-u)^{-1/r}-1}{k} \right\} \biggr)\\
&&s(u)=\left(\log r/u\right)^{-2} \frac{\varepsilon(u) d(u)}{\left(1 +\frac{\log k + \log d(u)}{\log r/u}\right)^2 + (\log r/u)^{-2}},
\end{eqnarray*}

\Bin with $\varepsilon(u)=(1-u)^{-(r+1)/r}$ , $d(u)=(u/r)/\{(1-u)^{-1/r}-1\}$, $arctan(\circ)$ is the inverse function of the tangent function $\tan(\circ)$

\begin{eqnarray*}
(Burr \ VI)&:& F^{-1}(1-u)= arcsinh \bigg( -\log\left(\frac {(1-u)^{-1/r}-1}{k}\right) \biggr)\\
&&s(u)=\left(\log r/u\right)^{-1} \frac{\varepsilon(u) d(u)}{\left(1 +\left(\frac{\log k +\log d(u)}{\log r/u}\right)^2 +(\log r/u)^{-2}\right)^{1/2}},
\end{eqnarray*}

\Bin with $\varepsilon(u)=(1-u)^{-(r+1)/r}$ , $d(u)=(u/r)/\{(1-u)^{-1/r}-1\}$, $arcsinh(\circ)$ is the inverse function of the hyperbolic sine function $\sinh(\circ)$

\begin{eqnarray*}
(Burr \ X)&:& F^{-1}(1-u)= \biggr(-\log\left((1-u)^{-1/r}-1\right)\biggr)^{1/2}\\
&&s(u)=\frac{1}{2}\left(\log r/u\right)^{-1/2} \frac{\varepsilon(u) d(u)}{\left(1 + \frac{\log d(u)}{\log r/u}\right)},
\end{eqnarray*}

\Bin with $\varepsilon(u)=(1-u)^{-(r+1)/r}$ , $d(u)=(u/r)/\{(1-u)^{-1/r}-1\}$.\\

\begin{eqnarray*}
(Burr \ Xa)&:& F^{-1}(1-u)= \biggr(-\log\left(1-(1-u)^{1/r}\right)\biggr)^{1/2}\\
&&s(u)=\frac{1}{2}\left(\log r/u\right)^{-1/2} \frac{\varepsilon(u) d(u)}{\left(1 + \frac{\log d(u)}{\log r/u}\right)},
\end{eqnarray*}

\Bin with $\varepsilon(u)=(1-u)^{(1-r)/r}$ , $d(u)=(u/r)/\{1-(1-u)^{1/r}\}$.\\

\Ni So, for all four cases, we have $s(\circ)$ is slowly varying at zero and hence by Representation \eqref{portal.rdgr} in Proposition \ref{portal.rd}, we have that
$F \in D(G_0)$. Also, $s(u)\rightarrow 0$ as $u\rightarrow 0$.\\

\Ni Now, we are going to use Theorem \ref{extensionBelow} to all four cases. We give the details for one case, the first for example ($Burr \ V $). We remark that
$\varepsilon(u)=1 + O(u)$ and $d(u)= O(u)$. For $min(v_n,V_n)\leq u \leq max(v_n,V_n)$, we have $A_n=\min(n,S_{(n)})\leq \log (1/u) \leq \max(n,S_{(n)})=B_n$. So uniformly in $(u,v) \in [min(v_n,V_n),  max(v_n,V_n)]^2$, we have for any $\eta>1$
$$
\varepsilon(u)=1 + O_{\mathbb{P}}(d_n(\eta)),  \ \ d(u)= O_{\mathbb{P}}(d_n(\eta)), \  \log u=O_{\mathbb{P}}(n),
$$

\Bin which leads to

$$
\frac{s(u)}{s(v)}=\left(\frac{\log (1/u) + \log r}{\log (1/v) + \log r}\right)^{-1/2} (1+O_{\mathbb{P}}(f_n)),
$$

\Bin for $f_n=n^{-2}$, since the rates $O_{\mathbb{P}}(n^{-1/2})$ are much lower that those of $O_{\mathbb{P}}(d_n(\eta))$. Since 

$$
A_n=min(n,S_{(n)})\leq \log (1/u), \log (1/v) \leq max(n,S_{(n)})=B_n,
$$

\Bin we have
$$
C_n=\left(\frac{\log A_n -\log r}{\log B_n -\log r}\right)^{-1/2} \leq \left(\frac{\log u -\log r}{\log v -\log r}\right)^{-1/2} \leq \left(\frac{\log B_n -\log r}{\log A_n -\log r}\right)^{-1/2}=D_n .
$$

\Bin It straightforward to show that $(C_n-1)$ and $(D_n-1)$ are both $O_{\mathbb{P}}(n^{-1/2})$. So the condition of Theorem \ref{extensionBelow} holds with 
$d_n=n^{-1/2}$. By handling the rates appropriately, we get (we recall that $c_n=(\log n)/\sqrt{n}$, $n\geq 1$): \label{al_records_Burr_V_alt}\\

\Ni \textbf{(V) Alternative form of the asymptotic law of record values\\from Burr V}. (with $f_n=n^{-2}$, $d_n=n^{-2}$)

\begin{eqnarray*}
\frac{(\log r + n)^{2} \biggr(X^{(n)}- \arctan \biggr(-\log\left\{\frac{(1-e^{-n})^{-1/r}-1}{k} \right\} \biggr)\biggr)}{\sqrt{n}}&=&S_{n}^{\ast} + O_{\mathbb{P}}(n^{-1/2})=W_n^\ast + O_{\mathbb{P}}(c_n)\notag\\
																											&\rightarrow &\mathcal{N}(0,1). \ \ \  \text{\textit{(VAlt)}}
\end{eqnarray*}

\Bin By using again the same techniques, we get the following forms.\\

\Ni \textbf{(VI) Alternative form of the asymptotic law of record values\\ from Burr VI}. (with $f_n=n^{-1/2}$, $d_n=n^{-1/2}$) \label{al_records_Burr_VI_alt}

\begin{eqnarray*}
\frac{(\log r + n)^{2} \biggr(X^{(n)}-arcsinh \biggr(-\log\left\{\frac{(1-e^{-n})^{-1/r}-1}{k} \right\} \biggr)\biggr)}{\sqrt{n}}&=&S_{n}^{\ast} + O_{\mathbb{P}}(n^{-1/2})=W_n^\ast + O_{\mathbb{P}}(c_n)\notag\\
																											&\rightarrow &\mathcal{N}(0,1). \ \ \  \text{\textit{(VIAlt)}}
\end{eqnarray*}

\Ni \textbf{(X) Alternative form of the asymptotic law of record values\\ from Burr X}. (with $f_n=n^{-1}$, $d_n=n^{-1/2}$) \label{al_records_Burr_X_alt}

\begin{eqnarray*}
\frac{(\log r + n)^{2} \biggr(X^{(n)}-\biggr(-\log\left((1-e^{-n})^{-1/r}-1\right)\biggr)^{1/2}\biggr)}{\sqrt{n}}&=&S_{n}^{\ast} + O_{\mathbb{P}}(n^{-1/2})=W_n^\ast + O_{\mathbb{P}}(c_n)\notag\\
&\rightarrow &\mathcal{N}(0,1). \ \ \  \text{\textit{(XAlt)}}
\end{eqnarray*}

\Ni \textbf{(Xa) Alternative form of the asymptotic law of record values\\ from Burr Xa}. (with $f_n=n^{-1/2}$, $d_n=n^{-1/2}$) \label{al_records_Burr_Xa_alt}

\begin{eqnarray*}
\frac{(\log r + n)^{2} \biggr(X^{(n)}-\biggr(-\log\left(1-(1-e^{n})^{-1/r}\right)\biggr)^{1/2}\biggr)}{\sqrt{n}}&=&S_{n}^{\ast} + O_{\mathbb{P}}(n^{-1/2})=W_n^\ast + O_{\mathbb{P}}(c_n)\notag\\
&\rightarrow &\mathcal{N}(0,1). \ \ \  \text{\textit{(XaAlt)}}
\end{eqnarray*}

\Bin \section*{Conclusions and Perspective}

\Ni This paper has its own interests by giving all the asymptotic laws of the upper records values that can be adapted to give the corresponding results for 
the lower records values. Actually, the results are at intersection of three major sub-disciplines of Statistics (Extreme value theory, records theory and asymptotic expansions) and the opportunity to summarize them has been seized. Most importantly, the study of the distributions of the so important Burr law paves the way of a handbook of similar results for as much as possible continuous distributions. This handbook already exists as a draft. The current paper will be main source of citations of it.\\

\Bin \textbf{Acknowledgment}.The authors wish to express their thanks to the anonymous reviewer for his valable and helpful comments. They also appreciated the efficient handling of the paper by the editor-in-chief.

\newpage

\begin{table}
	\centering
		\begin{tabular}{cccccccccccc}
	 \hline \hline
	 $n$			&25 		&50 		&75		&100 		&150 	&200 	&250 	&300 		&350 & 500 & 1000\\
	 
	$p(3)$		&0.45		&0.322	&0.26	&0.22		&0.18	&0.15	&0.14	&0.13		&0.11 &0.098&	0.069\\
	 \hline \hline
	 \end{tabular}
	\caption{Probabilities table of having less three records in a sample of size $n$}
	\label{Tab1}
\end{table}

\begin{table}
	\begin{centering}
		\begin{tabular}{|c|c|c|c|c|c|c|}
			\cline{2-7} \cline{3-7} \cline{4-7} \cline{5-7} \cline{6-7} \cline{7-7} 
			\multicolumn{1}{c|}{} & \multicolumn{1}{c}{$\gamma>0$} &  & \multicolumn{1}{c}{$\gamma<0$} &  & \multicolumn{1}{c}{$\gamma=0$} & \tabularnewline
			\hline 
			Burr & $II$ & $III$ & $I$ & $IV$ & $VI$ & $X$\tabularnewline
			$(parameters)$ & $(r=2)$ & $(r=2,k=3)$ & $ $ & $(r=2,c=3)$ & $(k=2,r=3)$ & $(k=2,r=3)$\tabularnewline
			\hline 
			$P_{0}(\%)$ & $3.78$ & $3.05$ & $2.28$ & $14.65$ & $3.54$ & $1.77$\tabularnewline
			\hline 
		\end{tabular}
		\par\end{centering}
	\caption{Empirical $p-$value of test record values normality ($nr=3$ records)}
	\label{Tab2}
\end{table}

\end{document}